\documentclass[openany,letterpaper,11pt]{amsart}

\usepackage{amsthm}
\usepackage{amsmath}
\usepackage{esint}
\usepackage{amssymb}
\usepackage{amsfonts}
\usepackage{graphicx}
\usepackage{ragged2e}
\usepackage{xcolor}
\usepackage{caption}
\usepackage{subcaption}

\usepackage[left=2.5cm,right=2.5cm,top=2cm,bottom=2cm]{geometry}
\theoremstyle{plain}

\newtheorem{proposition}{Proposition}
\newtheorem{remark}{Remark}

\numberwithin{equation}{section}

\usepackage{scalerel}
\usepackage{graphicx}
\begin{document}

\title[An interface problem for the Dirac equation]{Resolving an interface problem for the Dirac equation by using the unified transform method}
\author[C. A. Garc\'{i}a-Bibiano]{C. A. Garc\'{i}a-Bibiano}
\address{Facultad de Ciencias Qu\'{i}mico-Biol\'{o}gicas, Universidad Aut\'{o}noma de Guerrero, Av. L\'{a}zaro
C\'{a}rdenas S/N, Ciudad Universitaria, C.P. 39086, Chilpancingo de los Bravo, Guerrero, M\'{e}xico}
\email{carlosbibiano@uagro.mx} 

\begin{abstract}
\quad We use the Unified Transform Method (UTM) for the vector case to resolve an interface problem for the Dirac equation on two semi-infinite domains and two finite domains in the massless and massive cases, respectively. The UTM for the vector case is a variation of the UTM for the scalar case. The solutions obtained for an interface problem on two semi-infinite domains and two finite domains, respectively, in the massive case are convergent explicit integral representations.
\end{abstract}
\maketitle

\section{Introduction}
\label{section_1} 

\quad The interface problems in Partial Differential Equations (PDEs) are initial value and/or boundary problems for which the solution of the equation defined in a domain prescribes boundary conditions for equations in adjacent domains. In applications, interface problems (also known as transmission problems) frequently arise in heat transfer, quantum mechanics, and mathematical biology. Examples include heat flowing through a composite bar \cite{2,12}, the time-dependent linear Schrödinger equation with a piecewise constant potential \cite{7,14,17}, and shock waves as a method of healing broken bones \cite{6}. Finding solutions for equations that model water waves can also be considered an interface problem, in which the interface is composed of the contact region between air and water \cite{13,19}. In applications, many of the interface conditions are established to model a certain physical phenomenon that follows from conservation laws. To derive the boundary conditions that must be imposed at the interface, we consider differential equations that are valid on one side or the other of the interface. The integral form of the equation or important integral relations is often the best way to infer necessary conditions on unknown functions at the interface. In many physical systems modeled by differential equations, the system parameters may not be continuous functions of time and space. These discontinuities are often reflected in interface conditions for the differential equations; therefore, interface problems can be formulated in different PDEs. \\
The UTM or Fokas Method presents a new approach to solving Boundary Value Problems (BVPs) for nonlinear integrable differential equations \cite{8,10}, and it may even represent a more satisfactory and direct alternative for solving linear PDEs with constant coefficients \cite{5,10}. In the latter context, the method has provided new results for one-dimensional problems involving more than two spatial derivatives \cite{9,10}, for elliptical problems \cite{1,10}, and for interface problems \cite{4,Klein_Gordon,15,18}, among other references. \\
There are few studies on how the UTM is applied to systems of linear equations or to higher-order scalar equations that can be rewritten in the form of linear systems. The UTM for the vector case solves the Klein-Gordon equation, defined in $\mathbb{R}^{+}$, with initial conditions and Dirichlet boundary conditions (see \cite{3}). Another form, using Lax pairs, in which this initial value and/or boundary problem is solved, can be consulted in \cite{16}. On the other hand, using UTM for the vector case solves a BVP for the Dirac equation in $1$-dimension, see \cite{11}. \\
The Dirac equation can be solved using classical methods such as the Fourier series and the method of characteristics. In this paper, we employ the UTM, which is also known as Fokas's Method in the literature. The UTM is more powerful than the aforementioned classical methods because it is more applicable. We note that BVPs with the Schrödinger equation have been studied using the UTM, see \cite{Anne, Fokas1} among others. Recently, the method was applied to a general class of PDEs systems \cite{3}. \\
\quad In this paper, we study a family of interface problems for the Dirac equation in one spatial dimension,
\begin{equation}
\begin{cases}
\begin{split}
& i \partial_{t} \Psi_{1}^{(j)}(x,t)=-i \partial_{x} \Psi_{1}^{(j)}(x,t)+ m_{j} \Psi_{2}^{(j)}(x,t), && \ x \in \mathbb{R}^{-}, \textup{ if } j=1; \\
& i \partial_{t} \Psi_{2}^{(j)}(x,t)=i \partial_{x} \Psi_{2}^{(j)}(x,t) + m_{j} \Psi_{1}^{(j)}(x,t), &&  \ x \in \mathbb{R}^{+}, \textup{ if } j=2,
\label{ec2}
\end{split} 
\end{cases}
\end{equation}
where $t > 0$, $m_{j} \geq 0$, and $\displaystyle {\partial_{t} \Psi_{l}^{(j)}(x,t)}\equiv \frac{\partial \Psi_{l}^{(j)}(x,t)}{\partial t}$, for $l,j=1,2$; analogously for the other partial derivatives. In particular, we study the interface problems for this system on two semi-infinite domains $(- \infty , 0) \times (0, \infty )$ and two finite domains $[-L,0] \times [0,L]$, $L>0$. We show that the UTM for the vector case can be applied to obtain explicit solutions determined by the initial and boundary conditions of the problem posed. \\

\noindent The structure of this article is as follows: Section \ref{section_2} considers an interface problem for the Dirac equation on two semi-infinite domains. In Subsections \ref{section_2_1} and \ref{section_2_2}, we study the massless and massive systems, respectively. Finally, Section \ref{section_3} considers an interface problem on two finite domains. Again, in Subsections \ref{section_3_1} and \ref{section_3_2}, we study the massless and massive cases.

\section{Solution of an interface problem for the Dirac equation on two semi-infinite domains}
\label{section_2} 

\quad In this Section, we apply the UTM for the vector case to solve an interface problem for the Dirac equation on two semi-infinite domains, $\mathbb{R^{-}}:=\left\lbrace x \in \mathbb{R}:x<0 \right\rbrace$ and $ \mathbb{R^{+}}:=\left\lbrace x \in \mathbb{R}:x>0 \right\rbrace$. We seek four functions
\begin{equation}
\begin{array}{ccc}
\left.
\begin{split}
\Psi_{1}^{(1)}(x,t), \Psi_{2}^{(1)}(x,t), \ \ x  \in \mathbb{R}^{-}; \\
\Psi_{1}^{(2)}(x,t), \Psi_{2}^{(2)}(x,t), \ \  x \in \mathbb{R}^{+},  
\label{ec1}
\end{split} \right |
\ \ t > 0,
\end{array}
\end{equation}
\noindent functions that satisfy the equations
\begin{equation}
\begin{cases}
\begin{split}
& i \partial_t \Psi_1^{(j)}(x, t)=-i \partial_x \Psi_1^{(j)}(x, t)+ m_{j} \Psi_2^{(j)}(x, t), && \ x \in \mathbb{R}^{-}, \textup{ if } j=1; \\
& i \partial_t \Psi_2^{(j)}(x, t)=i \partial_x \Psi_2^{(j)}(x, t) + m_{j} \Psi_1^{(j)}(x, t), && \ x \in \mathbb{R}^{+}, \textup{ if } j=2,
\label{ec2}
\end{split} 
\end{cases}
\end{equation}
where $t > 0$, $m_{j} \in \mathbb{R}^{+}$, and $\displaystyle { \partial_t \Psi_{l}^{(j)}}(x,t)\equiv \frac{\partial \Psi_{l}^{(j)}(x,t)}{\partial t}$, for $j,l=1,2$; analogously for the other partial derivatives. \\
\noindent Suppose that we have the initial conditions
\begin{equation}
\begin{array}{ccc}
\Psi_{1}^{(1)}(x,0)=\Psi_{1,0}^{(1)}(x), \ \ x \in \mathbb{R}^{-};\\
\Psi_{1}^{(2)}(x,0)=\Psi_{1,0}^{(2)}(x), \ \ x \in \mathbb{R}^{+};\\
\Psi_{2}^{(1)}(x,0)=\Psi_{2,0}^{(1)}(x), \ \ x \in \mathbb{R}^{-}; \\
\Psi_{2}^{(2)}(x,0)=\Psi_{2,0}^{(2)}(x), \ \ x \in \mathbb{R}^{+},
\label{ec3}
\end{array}
\end{equation}
the asymptotic conditions
\begin{equation}
\begin{array}{ccc}
\left.
\begin{split}
\lim_{x \to -\infty} \Psi_{j}^{(1)}(x,t)= \lim_{x \to -\infty} \partial_x \Psi_{j}^{(1)}(x,t) = \lim_{x \to -\infty} \partial_t \Psi_{j}^{(1)}(x,t) = 0, \ j=1,2; \\
\lim_{x \to \infty} \Psi_{j}^{(2)}(x,t)= \lim_{x \to \infty} \partial_x \Psi_{j}^{(2)}(x,t) = \lim_{x \to \infty} \partial_t \Psi_{j}^{(2)}(x,t) = 0, \ j=1,2,  
\label{ec4}
\end{split} \right |
\ \ t \geq 0,
\end{array}
\end{equation}  
and the continuity conditions in the interface $x=0$,
\begin{equation}
\begin{array}{ccc}
\left.
\begin{split}
\Psi_{1}^{(1)}(0,t)=\Psi_{1}^{(2)}(0,t); \\
\Psi_{2}^{(1)}(0,t)=\Psi_{2}^{(2)}(0,t),  
\label{ec5}
\end{split} \right |
\ \ t > 0.
\end{array}
\end{equation}
To apply the Unified Transform Method to linear systems of partial differential equations, we follow the approach given in \cite{3} and rewrite \eqref{ec2} in the form,
\begin{equation}
\partial_{t}Q^{(j)}(x,t) + \Lambda^{(j)}(-i\partial_{x})Q^{(j)}(x.t)=0, \ \ j=1,2.
\label{vector_case}
\end{equation}
Here $Q^{(j)}(x,t)$ is an $N$-dimensional vector and $\Lambda^{(j)}(-i\partial_{x})$ is a $N \times N$ matrix-valued polynomial; $N = 2$ for the one-dimensional Dirac system. We now rewrite \eqref{ec2} in matrix form as,
\begin{equation}
\partial_t \displaystyle{\binom{\Psi_1^{(j)}(x,t)}{\Psi_2^{(j)}(x,t)}} + \left(\begin{array}{cc}
\partial_x & m_{j}i \\
m_{j}i & -\partial_x
\end{array}\right)\binom{\Psi_1^{(j)}(x,t)}{\Psi_2^{(j)}(x,t)}=0, \ j=1,2.
\label{matrix_form}
\end{equation}
Thus, from \eqref{vector_case}, we have that
\begin{equation}
Q^{(j)}(x,t) = \displaystyle{\binom{\Psi_1^{(j)}(x,t)}{\Psi_2^{(j)}(x,t)}}, \ \ \Lambda^{(j)}(k)=\left(\begin{array}{cc}
i k & m_{j}i \\
m_{j}i & -i k
\end{array}\right) =0, \ \ j=1,2.
\label{ecuacion1}
\end{equation}

\noindent The main distinction between the interface problem for the Dirac equation for the massive case $(m_{j} > 0, \ j=1,2)$ and the interface problem for the Dirac equation for the massless case $(m_{j} = 0, \ j=1,2)$ is the level of complexity of the EDPs system that must be resolved in each respective case.

\subsection{The massless system on two semi-infinite domains}
\label{section_2_1}

To solve the massless system, we set $m_{j}=0$, for $j=1,2$, in \eqref{matrix_form}. Thus, we have

\begin{equation*}
\partial_t \displaystyle{\binom{\Psi_1^{(j)}(x,t)}{\Psi_2^{(j)}(x,t)}} + \left(\begin{array}{cc}
\partial_x & 0 \\
0 & -\partial_x
\end{array}\right)\binom{\Psi_1^{(j)}(x,t)}{\Psi_2^{(j)}(x,t)}=0, \ j=1,2;
\end{equation*}

\noindent and

\begin{equation}
\Lambda(k)=\left(\begin{array}{cc}
i k & 0 \\
0 & -i k
\end{array}\right).
\label{ecuacion1}
\end{equation}
\\
\noindent We notice that the massless case is overdetermined because the boundary conditions of $\Psi_{2}^{(j)}$ for $j=1,2$, cannot be specified. One of the key aspects of the UTM is that it details how the boundary conditions are necessary for a well-posed problem. Applying the UTM, we show that the initial conditions of $\Psi_2^{(j)}$ dictate the boundary conditions of $\Psi_2^{(j)}$ for $j=1,2$. Thus, making it impossible to independently prescribe these boundary conditions adequately.

\noindent The first step to compute the solution to the massless system is to calculate the dispersion relations. Computing the dispersion relations for the vector case is analogous to applying the UTM over scalar PDEs, see \cite{3,5}. However, unlike the scalar case, where the dispersion relation is computed by solving the adjoint problem or by substituting $e^{i k x-\Omega(k)t}$ into the PDE following the convention given in \cite{5}, we define
\begin{equation}
Q^{(j)}(x,t)=\binom{\Psi_1^{(j)}(x,t)}{\Psi_2^{(j)}(x,t)} e^{i k x-\Omega(k) t}, \ j=1,2.
\label{Q_j}
\end{equation}

\noindent Then, in the system equivalent to the Dirac equation using \eqref{Q_j} to solving for $\Omega(k)$, we find $\Omega(k)$ such that
\begin{equation}
\operatorname{det}(\Lambda(k)-\Omega(k) I)=0 .
\label{ecuacion2}
\end{equation}

\noindent Computing the determinant, and solving for $\Omega(k)$, we determine branches to the massless system,
\begin{equation}
\Omega_{1,2}(k)= \pm ik.
\label{ecuacion3}
\end{equation}

\noindent Thus, the branches of $\Omega(k)$ are $\Omega_{1}(k)=ik$ and $\Omega_{2}(k)=-ik$. The next step in the vector case is to determine the local relation. In this case, we must rewrite (\ref{vector_case}) in divergence form (see the appendix of \cite{BVP}),
\begin{equation}
\left(e^{-i k x I+\Omega(k) t} A(k) Q^{(j)}(x,t)\right)_t-\left(e^{-ikx I+\Omega(k) t} A(k) X(x,t,k) Q^{(j)}(x,t)\right)_x=0, \ j=1,2;
\label{ecuacion4}
\end{equation}
where $\Omega(k)=\operatorname{diag}\left(\Omega_{1}(k), \ldots, \Omega_{N}(k)\right)$ is a diagonal matrix of the different branches of the system and $A(k)$ diagonalizes $\Lambda(k)$, that is,
$$
\Lambda(k)=A^{-1}(k) \Omega(k) A(k).
$$

\noindent Following \cite{3}, the matrix $X(x,t,k)$ is defined by
\begin{equation}
X(x,t,k)=\frac{i}{k+i \partial_x}\left(\Lambda(k)-\Lambda\left(-i \partial_x\right)\right).
\label{ecuacion5}
\end{equation}

\noindent For the massless Dirac equation, since $\Psi_{1,2}^{(j)}(k)=\Omega_{1,2}(k), \ j=1,2;$ and $A(k)=I$, we have
$$
X(x,t,k)=\frac{i}{k+i \partial_x}\left(\begin{array}{cc}
i k-\partial_x & 0 \\
0 & -i k+\partial_x
\end{array}\right)=\left(\begin{array}{cc}
-1 & 0 \\
0 & 1
\end{array}\right).
$$

\noindent We note that the simple form of $X(x,t,k)$ comes from the fact that when $m_{j}=0, \ j=1,2;$ the system is decoupled. This is not the case for the massive system. From (\ref{ecuacion4}), we now compute the local relations for the massless Dirac equation,
\begin{eqnarray}
\begin{split}
& \left(e^{-i k x+\Omega_1(k)t} \Psi^{(j)}_1(x,t)\right)_t-\left(e^{-i k x+\Omega_1(k)t}\left(-\Psi^{(j)}_1(x,t)\right)\right)_x=0, \  j=1,2;\\
& \left(e^{-i k x+\Omega_2(k)t} \Psi_2^{(j)}(x,t)\right)_t-\left(e^{-i k x+\Omega_2(k)t} \Psi_2^{(j)}(x,t)\right)_x=0, \ j=1,2.
\label{ecuacion6}
\end{split}
\end{eqnarray}

\noindent The next step is to compute the global relations from each of the local relations. In \cite{3}, a generalized expression was given to compute the global relations; however, the unknown variables in this problem are sufficiently decoupled so that the scalar process can be used to determine the global relations for (\ref{ecuacion6}). Computing the global relations is done by integrating each local relation over the corresponding domain, that is, we integrate \eqref{ecuacion6} on $C = (- \infty , 0] \times (0,T], \ T>0$, when $j=1$ and we integrate \eqref{ecuacion6} on $D = [0, \infty ) \times (0,T], \ T>0$ when $j=2$ with Green's theorem, respectively.
\begin{figure}[h]
\centering
\includegraphics[height=4.5cm,width=8cm]{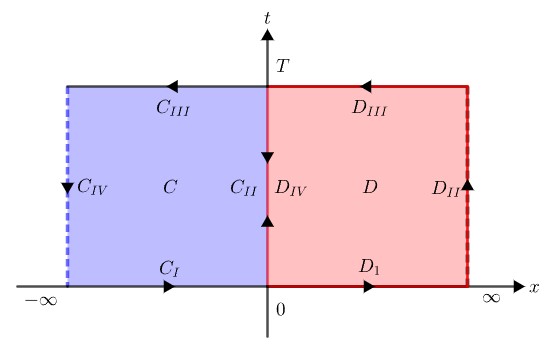}
\caption{Integration domains for the interface problem posed.}
\label{Figura_1}
\end{figure}

\noindent First, we compute the global relation for $\Psi_1^{(1)}$. Integrate the first equation in \eqref{ecuacion6} over $C$,
$$
\iint_{C}\left[e^{-ikx+\Omega_1(k)t} \Psi_1^{(1)}(x,t)\right]_t-\left[e^{-ikx+\Omega_1(k)t}\left(-\Psi_1^{(1)}(x,t)\right)\right]_x dtdx=0.
$$

\noindent Applying Green's Theorem, we move the integration to the boundary of $C$, 
$$
\int_{\partial C}\left[e^{-ikx+\Omega_1(k)t} \Psi_1^{(1)}(x,t)\right] dx+\left[e^{-ikx+\Omega_1(k)t}\left(-\Psi_1^{(1)}(x,t)\right)\right] dt=0.
$$

\noindent Parameterizing the borders of the domain and integrating the respective line integrals leads to the global relation. The form of the global relation is given in terms of the Fourier transform of the unknown variable and its initial conditions. In this paper, we use the following definitions for the Fourier transform, the inverse Fourier transform, and other Fourier transforms of important functions for the interface problem posed
$$
\begin{gathered}
\mathcal{F}\{f(x,t)\}=\hat{f}(k,t)=\int_{-\infty}^{\infty} e^{-ikx} f(x,t) dx, \textup{ for } \ x \in \mathbb{R}, \ t>0, \ k \in \mathbb{C}; \\
\mathcal{F}^{-1}\{\hat{f}(k,t)\}=f(x,t)=\frac{1}{2 \pi} \int_{-\infty}^{\infty} e^{ikx} \hat{f}(k,t) dk, \textup{ for } \ x \in \mathbb{R}, \ t>0; \\
\displaystyle \hat{f}^{(1)}(k,t)=\int \limits_{-\infty}^0 e^{-ikx}f^{(1)}(x,t)dx, \textup{ for }  x<0, \ t>0, \ \textup{Im}(k)>0; \\
\displaystyle \hat{f}^{(2)}(k,t)=\int \limits_{0}^{\infty} e^{-ikx}f^{(2)}(x,t)dx, \textup { for }  x>0, \ t>0, \ \textup{Im}(k)<0.
\end{gathered}
$$

\noindent \noindent The contribution of the second line integral is zero because $\Psi_1^{(1)}(x,t) \rightarrow 0$ sufficiently fast as $x \rightarrow - \infty$ for all $t=T>0$ and the line integral is over a finite interval since $T$ is finite. Thus, integrating the first equation in \eqref{ecuacion6} over the contributions $C_{j}$, $j \in \{ I,II,III,IV \}$, that make up $\partial C$, we have the following: \\
$$
\hat{\Psi}_1^{(1)}(k,0)-e^{\Omega_1(k)T} \hat{\Psi}_1^{(1)}(k,T)-g_0\left(\Omega_1(k), \Psi_{1}^{(1)}(0,T), T\right)=0,
$$
where we define the variables $g_0\left(\Omega_j(k), \Psi_{j}^{(1)}(0,t), t\right)$, $j \in \{1,2\}$, which depend on the boundary conditions as
\begin{equation}
g_0\left(\Omega_j(k), \Psi_{j}^{(1)}(0,t), t\right) = \int_0^t e^{\Omega_j(k) s} \Psi_j^{(1)}(0,s)ds=\int_0^t e^{\Omega_j(k) s} \Psi_j^{(2)}(0,s)ds.
\label{ecuacion7}
\end{equation}

\noindent We also note that the negative half-line definition of the Fourier transform is utilized in the calculation because the values of $\Psi_1^{(1)}(x,t)$ and $\Psi_2^{(1)}(x,t)$ are assumed to be zero for all $(x,t)$ outside the domain $C$. Following the same procedure for $\Psi_2^{(1)}(x,t)$ leads to the second global relation,
$$
\hat{\Psi}_{2}^{(1)}(k,0)+e^{\Omega_2(k) T} \hat{\Psi}_2^{(1)}(k,T)-g_0\left(\Omega_2(k), \Psi_2^{(1)}(0,T), T\right)=0
$$

\noindent We can replace $T$ in global relations with any $t \in (0,T]$ because both equations are valid for all feasible values of $t$. Thus, the global relations for the massless case where $j=1$ are
\begin{eqnarray}
\begin{split}
& \hat{\Psi}_{1}^{(1)}(k,0)-e^{\Omega_1(k)t} \hat{\Psi}_1^{(1)}(k,t)-g_0\left(\Omega_1(k), \Psi_{1}^{(1)}(0,t), t\right)=0, \\
& \hat{\Psi}_2^{(1)}(k,0)+e^{\Omega_2(k) t} \hat{\Psi}_2^{(1)}(k,t)-g_0\left(\Omega_2(k), \Psi_{2}^{(1)}(0,t), t\right)=0.
\label{ecuacion8_1}
\end{split}
\end{eqnarray}

\noindent For the other part, integrating the local relation \eqref{ecuacion6} in the domain $D = (0, \infty) \times (0,T), \ T > 0$, when $j=2$, and using Green's Theorem, that is, integrating the equations \eqref{ecuacion6} over the contributions $D_{j}$, $j \in \{ I,II,III,IV \}$, when $j=2$, respectively, we obtain
\begin{eqnarray}
\begin{split}
& \hat{\Psi}_{1}^{(2)}(k,0)-e^{\Omega_1(k)t} \hat{\Psi}_{1}^{(2)}(k,t)+g_0\left(\Omega_1(k), \Psi_{1}^{(2)}(0,t), t\right)=0; \\
& \hat{\Psi}_{2}^{(2)}(k,0)-e^{\Omega_2(k) t} \hat{\Psi}_2^{(2)}(k,t)-g_0\left(\Omega_2(k), \Psi_{2}^{(2)}(0,t), t\right)=0,
\label{ecuacion8_2}
\end{split}
\end{eqnarray}
where we define the variables $g_0\left(\Omega_{j}(k), \Psi_{j}^{(2)}(0,t), t\right)$, $j=1,2;$, which depend on the boundary conditions as
\begin{equation}
g_{0}\left(\Omega_{j}(k), \Psi_{j}^{(2)}(0,t), t\right) = \int_{0}^{t} e^{\Omega_j(k) s} \Psi_j^{(2)}(0,s)ds=\int_0^t e^{\Omega_{j}(k) s} \Psi_j^{(1)}(0,s)ds.
\label{ecuacion7}
\end{equation}

\noindent The next step is to utilize the inverse Fourier transform to obtain expressions for $\Psi_1^{(1)}(x,t)$ and $\Psi_2^{(1)}(x, t)$ from \eqref{ecuacion8_1}. Multiplying the first equation in \eqref{ecuacion8_1} by $e^{-\Omega_1(k)t}$ and the second by $e^{-\Omega_2(k)t}$, then algebraically solving for the Fourier transforms yields
$$
\begin{aligned}
& \hat{\Psi}_1^{(1)}(k,t)=e^{-\Omega_1(k)t}\left(\hat{\Psi}_{1,0}^{(1)}(k)-g_0\left(\Omega_1(k), \Psi_{1}^{(1)}(0,t), t\right)\right), \\
& \hat{\Psi}_2^{(1)}(k,t)=-e^{-\Omega_2(k)t}\left(\hat{\Psi}_{2,0}^{(1)}(k)-g_0\left(\Omega_2(k), \Psi_{2}^{(1)}(0,t), t\right)\right).
\end{aligned}
$$

\noindent Taking the inverse Fourier transforms of the above expressions and substituting in the dispersion relations \eqref{ecuacion3}, we have
\begin{eqnarray}
\begin{split}
& \Psi_{1}^{(1)}(x,t)=\frac{1}{2 \pi} \int_{-\infty}^{\infty} e^{i k(x-t)}\left(\hat{\Psi}_{1,0}^{(1)}(k)-g_0\left(\Omega_1(k), \Psi_{1}^{(1)}(0,t), t\right)\right) dk, \\
& \Psi_{2}^{(1)}(x,t)=-\frac{1}{2 \pi} \int_{-\infty}^{\infty} e^{i k(x+t)}\left(\hat{\Psi}_{2,0}^{(1)}(k)-g_0\left(\Omega_2(k), \Psi_{2}^{(1)}(0,t), t\right)\right) dk.
\label{ecuacion9}
\end{split}
\end{eqnarray}

\noindent Proceeding similarly to \eqref{ecuacion8_2}, we found $\Psi_1^{(2)}(x,t)$ and $\Psi_2^{(2)}(x,t)$, respectively. Thus, we have
\begin{eqnarray}
\begin{split}
& \Psi_{1}^{(2)}(x,t)=\frac{1}{2 \pi} \int_{-\infty}^{\infty} e^{i k(x-t)}\left(\hat{\Psi}_{1,0}^{(2)}(k)+g_0\left(\Omega_{1}(k), \Psi_{1}^{(2)}(0,t),t\right)\right) dk, \\
& \Psi_{2}^{(2)}(x,t)=\frac{1}{2 \pi} \int_{-\infty}^{\infty} e^{i k(x+t)}\left(\hat{\Psi}_{2,0}^{(2)}(k)-g_0\left(\Omega_{2}(k), \Psi_{2}^{(2)}(0,t),t\right)\right) dk.
\label{ecuacion12}
\end{split}
\end{eqnarray}

\noindent Now, \eqref{ecuacion9} is not a solution because the expressions are overdetermined. So, our next step is to remove the unnecessary boundary condition of $\Psi_2^{(1)}$ and simplify \eqref{ecuacion9} to arrive at the solution. Normally, when applying UTM, one accomplishes this by deforming the paths of integration, respectively. In this case, the dispersion relations have zero non-trivial discrete symmetries, so we use analytic properties to simplify \eqref{ecuacion9}. Expanding \eqref{ecuacion9} with the definitions for $g_0\left(\Omega_j(k), \Psi_{j}^{(1)}(0,t), t\right)$ given in \eqref{ecuacion7},
\begin{eqnarray}
\begin{split}
& \Psi_1^{(1)}(x,t)=\frac{1}{2 \pi} \int_{-\infty}^{\infty} e^{ik(x-t)} \hat{\Psi}_{1,0}^{(1)}(k) dk - \frac{1}{2 \pi} \int_{-\infty}^{\infty} \int_0^t e^{ik(x-t+s)} \Psi_1^{(1)}(0,s) ds dk, \\
& \Psi_2^{(1)}(x,t)= - \frac{1}{2 \pi} \int_{-\infty}^{\infty} e^{ik(x+t)} \hat{\Psi}_{2,0}^{(1)}(k) dk + \frac{1}{2 \pi} \int_{-\infty}^{\infty} \int_0^t e^{ik(x+t-s)} \Psi_2^{(1)}(0,s) ds dk.
\label{ecuacion10}
\end{split}
\end{eqnarray}
\begin{figure}[h]
\centering
\includegraphics[height=5cm,width=8.75cm]{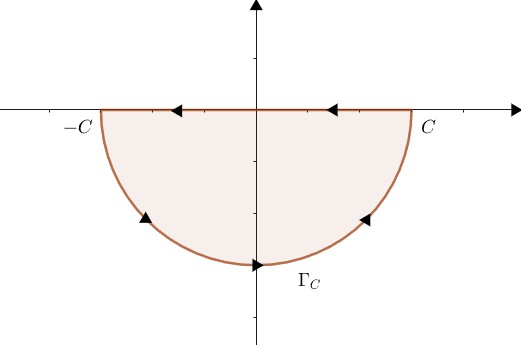}
\caption{Contour used to remove the boundary conditions for $\Psi_{1}^{(1)}(x,t)$ in \eqref{ecuacion10}}.
\label{Figura_2}
\end{figure}

\noindent In the first equation of \eqref{ecuacion10}, we see the integrand for the second integral over the real line is analytic. Since $x -t < 0$, then $x-t+s<0$, we integrate 
$$
\int_0^t e^{ik(x-t+s)} \Psi_1^{(1)}(0,s) ds
$$
along the contour $\Gamma_{C}=[-C,C] \cup \textup{Arc}_{C}$, where $\textup{Arc}_{C}$ is the circular arc of radius $C$ centered at the origin in $\mathbb{C}^{-}$, by Cauchy's Theorem, we have
$$
\left[ \int_{-C}^{C} + \int_{\textup{Arc}_{C}} \right] \left( \int_0^t e^{ik(x-t+s)} \Psi_1^{(1)}(0,s)ds \right)dk = 0.
$$
Taking the limit as $C \rightarrow \infty$ and applying Jordan's Lemma, due to the exponential decay, the integral along $\Gamma_{C}$ goes to $0$; therefore,
$$
\int_{-\infty}^{\infty} \int_0^t e^{ik(x-t+s)} \Psi_1^{(1)}(0,s)ds dk = 0,
$$
and using Fourier inversion, we obtain 
$$
\Psi_{1}^{(1)}(x,t) = \frac{1}{2 \pi} \int_{-\infty}^{\infty} e^{ik(x-t)} \Psi_{1,0}^{(1)}(k)dk = \Psi_{1,0}^{(1)}(x-t).
$$

\noindent To simplify $\Psi_{2}^{(1)}(x,t)$ in \eqref{ecuacion10}, we work in cases. \\

\noindent ${\bf{Case \ 1}}$: If $-x > t$, then $x + t - s < 0$ and using the previous arguments, we have
$$
\Psi_{2}^{(1)}(x,t)= - \frac{1}{2 \pi} \int_{-\infty}^{\infty} e^{i k(x+t-s)} \Psi_{2,0}^{(1)}(k) dk = - \Psi_{2,0}^{(1)}(x+t). \\
$$

\noindent ${\bf{Case \ 2}}$: If $-x < t$, then $x + t - s > 0$, and we cannot use contour $\Gamma_{c}$ to eliminate the second integral that defines $\Psi_{2}^{(1)}(x,t)$. However, using the distributional definition of the delta distribution, we have
$$
\int_{-\infty}^{\infty} \int_{0}^t e^{ik(x+t-s)} \Psi_{2}^{(1)}(0,s)ds dk = \Psi_{2}^{(1)}(0,x+t).
$$
So, for $-x<t$,
$\Psi_{2}^{(1)}(x,t) = \Psi_{2}^{(1)}(0,x+t)$, \\

\noindent by the definition of the inverse Fourier transform and $\Psi_{2,0}^{(1)}(x+t)=0$, since $x+t>0$. \\

\noindent In this work, we assume that all boundary and initial conditions are compatible whenever they interface. Therefore, in this problem, when $x=t$, $\Psi_{1}^{(j)}(x-t)=\Psi_{1}^{(j)}(0,t-x)$, $j=1,2$.

\begin{remark}
We note that the delta distribution argument in {\bf{Case 2}} could have been utilized for {\bf{Case 1}} instead of using the contour $\Gamma_{C}$, however, this is only due to the simple dispersion relation \eqref{ecuacion3} for the massless Dirac equation. Since invoking the delta distribution does not generalize to determining the massive solution, we use the contour argument.
\end{remark}

\noindent Using an analogous procedure like the one given above to solve the interface problem for the massless Dirac equation on the negative half-line, and on the other hand, using similar arguments and subsequently working by cases, the expressions obtained for $\Psi_{1}^{(2)}(x,t)$ and $\Psi_{2}^{(2)}(x,t)$, to solve the interface problem for the massless Dirac equation on the positive half-line. We obtain the following result:

\begin{proposition}
The solution of the interface problem for the Dirac equation \eqref{ec1}-\eqref{ec5} for the massless case is given by  
\begin{eqnarray}
\begin{split}
& \Psi_{1}^{(1)}(x,t)= \Psi_{1,0}^{(1)}(x-t),
\label{ecuacion11}
\end{split}
\end{eqnarray}
\noindent and
\begin{eqnarray}
\begin{split}
& \Psi_{2}^{(1)}(x,t)= \begin{cases} -\Psi_{2,0}^{(1)}(x+t), & -x \geq t; \\
\Psi_{2}^{(1)}(0,x+t), & -x<t , \end{cases} \\
\label{ecuacion12}
\end{split}
\end{eqnarray}
for $x<0$,\ $t>0$. \\
\begin{eqnarray}
\begin{split}
& \Psi_{1}^{(2)}(x,t)= \begin{cases}\Psi_{1,0}^{(2)}(x-t), & x \geq t; \\
\Psi_{1}^{(1)}(0,t-x), & x<t , \end{cases} \\
\label{ecuacion13}
\end{split}
\end{eqnarray}
\noindent and
\begin{eqnarray}
\begin{split}
& \Psi_{2}^{(2)}(x,t)=\Psi_{2,0}^{(2)}(x+t),
\label{ecuacion14}
\end{split}
\end{eqnarray}
for $x>0$,\ $t>0$.
\end{proposition}

\begin{remark}
Notice, since $\Psi_{1,0}^{(1)}(x-t)= \Psi_{1}^{(1)}(0,t-x)$ and $\Psi_{1,0}^{(2)}(x-t)=\Psi_{1}^{(2)}(0,t-x)$ when $x=t$ the equality in \eqref{ecuacion11} and \eqref{ecuacion13}, respectively, could have been in either case. Now, the forms of these solutions were expected because the massless case in two semi-infinite domains \eqref{ec1}-\eqref{ec5} reduces to the transport equation for $\Psi_{1}^{(1)}(x,t)$, $\Psi_{2}^{(1)}(x,t)$, $\Psi_{1}^{(2)}(x,t)$, and $\Psi_{2}^{(2)}(x,t)$ in their respective domains of definition. The method of characteristics provides the solution to the transport equation as left- and right-moving waves, which agree with \eqref{ecuacion11} and \eqref{ecuacion13}, respectively.
\end{remark}

\subsection{The massive system on two semi-infinite domains}
\label{section_2_2}

Now, we solve the interface problem for the massive Dirac equation on two semi-infinite domains. For this, we set $m_{j}>0$, $j=1,2$ in \eqref{ec2} giving
\begin{eqnarray}
\left\{\begin{array}{l}
i \partial_t \Psi_1^{(j)}(x,t)=-i \partial_x \Psi_1^{(j)}(x,t ) +m_{j} \Psi_2^{(j)}(x,t); \\
i \partial_t \Psi_2^{(j)}(x,t)=i \partial_x \Psi_2^{(j)}(x,t) + m_{j} \Psi_1^{(j)}(x,t), \\
\end{array}\right.
\label{ec6}
\end{eqnarray}
\noindent where if $j=1$, \eqref{ec6} is defined for $x \in \mathbb{R}^{-}$ and $t > 0$, while if $j=2$, \eqref{ec6} is defined for $x \in \mathbb{R}^{+}$ and $t > 0$. To apply the UTM for the vector case, we follow the approach given in \cite{3} and rewrite (\ref{ec6}) in the form
\begin{equation}
\partial_t Q^{(j)}(x,t)+\Lambda^{(j)}\left(-i \partial_x\right) Q^{(j)}(x,t)=0, \ j=1,2.
\label{ec7}
\end{equation}

\noindent Here $Q^{(j)}(x,t)$ is an $N$-dimensional vector and $\Lambda^{(j)}\left(-i \partial_x\right)$ is a $N \times N$ matrix-valued polynomial; $N=2$ for the one-dimensional Dirac system. We rewrite (\ref{ec6}) in matrix form as
\begin{equation}
\partial_t \displaystyle{\binom{\Psi_1^{(j)}(x,t)}{\Psi_2^{(j)}(x,t)}} + \left(\begin{array}{cc}
\partial_x & i m_{j} \\
i m_{j} & -\partial_x
\end{array}\right)\binom{\Psi_1^{(j)}(x,t)}{\Psi_2^{(j)}(x,t)}=0, \ j=1,2. 
\label{ec8}
\end{equation}

\noindent Thus, from (\ref{ec7}), we have that
\begin{equation}
Q^{(j)}(x,t)=\binom{\Psi_1^{(j)}(x, t)}{\Psi_2^{(j)}(x, t)} \ \textup{ and } \ \Lambda^{(j)}(k)=\left(\begin{array}{cc}
i k & i m_{j} \\
i m_{j} & -i k
\end{array}\right), \ j=1,2.
\label{ec9}
\end{equation}

\noindent However, $m_{j}>0$ ensures the system is coupled, and $\Lambda^{(j)}(k)$ is not a diagonal matrix. For these reasons, the massive system is more complicated than the massless; therefore, we solved the simpler massless case to introduce the technique before solving the massive case. Thus, in this case (\ref{ec6}) is not overdetermined, but the UTM will show in the coming steps that the solutions to $\Psi_{1}^{(j)}$ and $\Psi_{2}^{(j)}$ do not depend on each other's boundary conditions. In other words, the solution for $\Psi_{1}^{(j)}$ depends on both initial conditions and its boundary conditions, but not the boundary condition of $\Psi_{2}^{(j)}$ and vice versa for $\Psi_2^{(j)}$.

\noindent First, we begin by computing the dispersion relation. The computation follows as before in the massless case. We find $\Omega^{(j)}(k)$ such that,
$$
\operatorname{det}(\Lambda^{(j)}(k)-\Omega^{(j)}(k) I)=0, \ j=1,2.
$$
$\Lambda^{(j)}(k)$ is no longer a diagonal matrix, adding complexity to the dispersion relation. Therefore, performing the calculations, we have that $\Omega_{1,2}^{(j)}(k) = \pm i \sqrt{k^{2}+m_{j}^{2}}$, that is, the branches for the massive system are:
$$
\Omega_{1}^{(j)}(k) = i \sqrt{k^2+m_{j}^2} \textup{ \ and \ } \Omega_{2}^{(j)}(k)= - i \sqrt{k^2+m_{j}^2}, \ j=1,2.
$$

\noindent Thus, the branches of $\Omega^{(j)}(k)$ $(\Omega^{(j)}_{1}(k) \ \textup{and} \ \Omega^{(j)}_{2}(k), \ j=1,2)$, are roots of a polynomial in $k$ of order $2$. Obtaining the local relation for the massive system still requires the system to be in divergence form, given by
$$
\left(e^{-i k x I+\Omega^{(j)}(k) t} A^{(j)}(k) Q^{(j)}(x,t)\right)_t-\left(e^{-i k x I+\Omega^{(j)}(k) t} A^{(j)}(k) X(x, t, k) Q^{(j)}(x,t)\right)_x=0, \ j=1,2. 
$$

\noindent The value of $\displaystyle X(x,t,k)=\frac{i}{k+i\partial_{x}}\Big ( \Lambda^{(j)}(k) - \Lambda^{(j)}(-i\partial_{x}) \Big )$, does not change from the massless case, but the diagonalization of $\Lambda^{(j)}(k)$ is different due to the different branches. For the massive case,
$$
A^{(j)}(k)=\left(\begin{array}{ll}
i m_{j} & \Omega_1^{(j)}(k)-i k \\
i m_{j} & \Omega_2^{(j)}(k)-i k
\end{array}\right), \ j=1,2, \textup{ and } X(x,t,k)=\left(\begin{array}{cc}
-1 & 0 \\
0 & 1
\end{array}\right).
$$

\noindent With these definitions, after simplifying, we obtain the local relations for $j = 1,2$ 
\begin{equation}
\begin{array}{ccc}
\begin{split}
& 0 = \left(e^{-i k x+\Omega_{j}^{(1)}(k) t}\left[(im_{1}) \Psi_1^{(1)}(x,t)+\left(\Omega_{j}^{(1)}(k)-i k\right) \Psi_2^{(1)}(x,t)\right]\right)_t \\
& -\left(e^{-i k x+\Omega_{j}^{(1)}(k) t}\left[(-im_{1}) \Psi_1^{(1)}(x,t)+\left(\Omega_{j}^{(1)}(k)-i k\right) \Psi_2^{(1)}(x,t)\right]\right)_x, \ x < 0, \ t > 0;
\label{RL1_Dirac}
\end{split}
\end{array}
\end{equation}
\begin{equation}
\begin{array}{ccc}
\begin{split}
& 0 = \left(e^{-i k x+\Omega_{j}^{(2)}(k) t}\left[(im_{2}) \Psi_1^{(2)}(x,t)+\left(\Omega_{j}^{(2)}(k)-i k\right) \Psi_2^{(2)}(x,t)\right]\right)_t \\
& -\left(e^{-i k x+\Omega_{j}^{(2)}(k) t}\left[(-im_{2}) \Psi_1^{(2)}(x,t)+\left(\Omega_{j}^{(2)}(k)-i k\right) \Psi_2^{(2)}(x,t)\right]\right)_x, \ x > 0, \ t > 0.
\label{RL2_Dirac}
\end{split}
\end{array}
\end{equation}
\\
\noindent The global relations for the massive case are determined in the same manner as before in the massless case. We integrate the local relations \eqref{RL1_Dirac} on the domain $C = (-\infty,0) \times (0,T), \ T > 0$,
\begin{eqnarray*}
0 = & \displaystyle \iint_C\left(e^{-i k x+\Omega_{j}^{(1)}(k)t}\left[(im_{1}) \Psi_1^{(1)}(x,t)+\left(\Omega_{j}^{(1)}(k)-ik\right) \Psi_2^{(1)}(x,t)\right]\right)_t \\
& - \left(e^{-i k x+\Omega_{j}^{(1)}(k) t}\left[(-i m_{1}) \Psi_{1}^{(1)}(x,t)+\left(\Omega_{j}^{(1)}(k)-i k\right) \Psi_2^{(1)}(x,t)\right]\right)_x d x dt, \ j=1,2;
\end{eqnarray*}
applying Green's Theorem to move the integration to the boundary, we have
\begin{eqnarray*}
0 = & \displaystyle \int_{\partial C}\left(e^{-i k x+\Omega_{j}^{(1)}(k) t}\left[(i m_{1}) \Psi_1^{(1)}(x,t)+\left(\Omega_{j}^{(1)}(k)-i k\right) \Psi_2^{(1)}(x,t) \right]\right) dx \\
& +\left(e^{-i k x+\Omega_{j}^{(1)}(k) t}\left[(-i m_{1}) \Psi_1^{(1)}(x,t)+\left(\Omega_{j}^{(1)}(k)-i k\right) \Psi_2^{(1)}(x,t) \right]\right) dt, \ j=1,2.
\end{eqnarray*}
Then, parameterizing $\partial C = C_{I} + C_{II} + C_{III} + C_{IV}$, as in the massless case, and finally, using the properties of $\Psi_1^{(1)}$ and $\Psi_2^{(1)}$ to compute the various line integrals, we perform the integration that allows us to write the two global relations. The global relations for the massive case on $C$ are,
\begin{eqnarray}
\begin{split}
0 = & \left(i k-\Omega_{j}^{(1)}(k)\right) e^{\Omega_{j}^{(1)}(k) t} \hat{\Psi}^{(1)}_2(k,t)+(-im_{1}) h_{0,1}^{(1)}\left(\Omega_{j}^{(1)}(k), t\right) +(-im_{1}) e^{\Omega_{j}^{(1)}(k) t} \hat{\Psi}^{(1)}_1(k,t) \\
& + (im_{1}) \hat{\Psi}^{(1)}_1(k,0) + \left(\Omega_{j}^{(1)}(k)-i k\right) \hat{\Psi}^{(1)}_2(k,0) +\left(\Omega_{j}^{(1)}(k) - ik\right) h_{0,2}^{(1)}\left(\Omega_{j}^{(1)}(k), t\right), 
\label{ec11_1}
\end{split}
\end{eqnarray}
for $j = 1,2;$ where $T$ has been replaced by $t$ and
\begin{equation}
h_{0,l}^{(1)}\left(\Omega_{j}^{(1)}(k),t\right) = \int_0^t e^{\Omega_{j}^{(1)}(k) s} \Psi^{(1)}_l(0,s) ds, \ l=1,2.
\label{ec12}
\end{equation}
On the other hand, integrating the local relation \eqref{RL2_Dirac} on the domain $D = (0, \infty) \times (0,T), \ T > 0$, later, using Green's Theorem and parameterizing $\partial D = D_{I} + D_{II} + D_{III} + D_{IV}$, we obtain
\begin{eqnarray}
\begin{split}
0 = & \left(i k-\Omega^{(2)}_{j}(k)\right) e^{\Omega^{(2)}_{j}(k) t} \hat{\Psi}^{(2)}_2(k,t)+(im_{j}) h^{(2)}_{0,1}\left(\Omega^{(2)}_{j}(k), t\right) +(-im_{2}) e^{\Omega^{(2)}_{j}(k) t} \hat{\Psi}^{(2)}_1(k,t) \\
& + (im_{2}) \hat{\Psi}^{(2)}_1(k,0) + \left(\Omega^{(2)}_{j}(k)-i k\right) \hat{\Psi}^{(2)}_2(k,0) +\left(i k-\Omega^{(2)}_{j}(k)\right) h^{(2)}_{0,2}\left(\Omega^{(2)}_{j}(k), t\right),
\label{ec11_2}
\end{split}
\end{eqnarray}
for $j = 1,2,$ where $T$ has been replaced by $t$ and
\begin{equation}
h^{(2)}_{0,l}\left(\Omega^{(2)}(k),t\right) = \int_0^t e^{\Omega^{(2)}(k) s} \Psi^{(2)}_l(0,s) ds, \ l=1,2.
\label{ec12}
\end{equation}

\noindent To apply the inverse Fourier transform, we take the two global relations and solve for $\hat{\Psi}^{(1)}_1(k,t)$ and $\hat{\Psi}^{(1)}_2(k,t)$. We multiply equation (\ref{ec11_1}), with $j=1$ and $j=2$, by $e^{-\Omega^{(1)}_1(k) t}$ and $e^{-\Omega^{(1)}_2(k) t}$, respectively, to obtain,
\begin{eqnarray}
\begin{split}
& (-im_{1}) \hat{\Psi}^{(1)}_1(k,t)+\left(i k-\Omega^{(1)}_1(k)\right) \hat{\Psi}^{(1)}_2(k,t)+C^{(1)}_1(k,t) e^{-\Omega^{(1)}_1(k)t}=0, \\
& (-im_{1}) \hat{\Psi}^{(1)}_1(k,t)+\left(i k-\Omega^{(1)}_2(k)\right) \hat{\Psi}^{(1)}_2(k,t)+C^{(1)}_2(k,t) e^{-\Omega^{(1)}_2 (k)t}=0,
\label{ec13}
\end{split}
\end{eqnarray}
\noindent where
\begin{eqnarray*}
C^{(1)}_{n}(k,t) = & (-im_{1}) h_{0,1}^{(1)}\left(\Omega_{n}^{(1)}(k),t\right)+\left(\Omega_{n}^{(1)}(k)-ik\right) h_{0,2}^{(1)}\left(\Omega_{n}^{(1)}(k),t\right) \\
& + (im_{1}) \hat{\Psi}_{1,0}^{(1)}(k)+\left(\Omega_{n}^{(1)}(k)-ik\right) \hat{\Psi}^{(1)}_{2,0}(k), \ n = 1,2.
\end{eqnarray*}

\noindent Similarly, proceed to \eqref{ec11_2} and then applying the inverse Fourier transform to find $\Psi_1^{(2)}(x,t)$ and $\Psi_2^{(2)}(x,t)$, respectively. We have
\begin{eqnarray}
\begin{split}
& (-im_{2}) \hat{\Psi}^{(2)}_1(k,t)+\left(i k-\Omega^{(2)}_1(k)\right) \hat{\Psi}^{(2)}_2(k,t)+C^{(2)}_1(k,t) e^{-\Omega^{(2)}_1(k)t}=0, \\
& (-im_{2}) \hat{\Psi}^{(2)}_1(k,t)+\left(i k-\Omega^{(2)}_2(k)\right) \hat{\Psi}^{(2)}_2(k,t)+C^{(2)}_2(k,t) e^{-\Omega^{(2)}_2 (k)t}=0,
\label{ec13}
\end{split}
\end{eqnarray}
\noindent where
\begin{eqnarray*}
C^{(2)}_{r}(k,t) = & (im_{2}) h^{(2)}_{0,1}\left(\Omega_{r}^{(2)}(k),t\right)+\left(ik-\Omega_{r}^{(2)}(k)\right) h^{(2)}_{0,2}\left(\Omega_{r}^{(2)}(k),t\right) \\
& + (im_{2}) \hat{\Psi}_{1,0}^{(2)}(k)+\left(\Omega_{r}^{(2)}(k)-ik\right) \hat{\Psi}^{(2)}_{2,0}(k), \ r=1,2.
\end{eqnarray*}

\noindent Subtracting the second equation in \eqref{ec13} from the first, later, solving for $\hat{\Psi}^{(1)}_2(k,t)$, and using the fact $\Omega^{(1)}_2(k)-\Omega^{(1)}_1(k) = 2 \Omega^{(1)}_2(k)$, yields 
\begin{equation}
\hat{\Psi}^{(1)}_2(k,t)=\frac{1}{2 \Omega^{(1)}_2(k)}\left(C^{(1)}_2(k,t) e^{-\Omega^{(1)}_2(k) t}-C^{(1)}_1(k,t) e^{-\Omega^{(1)}_1(k)t}\right).
\label{ec14}
\end{equation}

\noindent Adding both equations in \eqref{ec13}, substituting \eqref{ec14}, and solving for $\hat{\Psi}^{(1)}_1(k, t)$,
\begin{equation}
\hat{\Psi}^{(1)}_1(k,t)=\left(\frac{\Omega^{(1)}_2(k) i+k}{2 m_{1} \Omega^{(1)}_1(k)}\right) C^{(1)}_1(k,t) e^{-\Omega^{(1)}_1(k) t}+\left(\frac{\Omega^{(1)}_2(k)i-k}{2 m_{1} \Omega^{(1)}_1(k)}\right) C^{(1)}_2(k,t) e^{-\Omega^{(1)}_2(k)t}.
\label{ec15}
\end{equation}

\noindent Taking the inverse Fourier transforms in \eqref{ec14} and \eqref{ec15}, respectively, we obtain
\begin{eqnarray}
\begin{split}
\Psi^{(1)}_1(x,t) = & \frac{1}{4 m_{1} \pi} \int_{-\infty}^{\infty}\left(\frac{\Omega^{(1)}_2(k) i+k}{\Omega^{(1)}_1(k)}\right) C^{(1)}_1(k,t) e^{i k x-\Omega^{(1)}_1(k) t} dk \\
& + \frac{1}{4 m_{1} \pi} \int_{-\infty}^{\infty}\left(\frac{\Omega^{(1)}_2(k) i-k}{\Omega^{(1)}_1(k)}\right) C^{(1)}_2(k,t) e^{i k x-\Omega^{(1)}_2(k)t} dk,
\label{ec16_1}
\end{split}
\end{eqnarray}
\begin{eqnarray}
\begin{split}
\Psi^{(1)}_2(x,t)=\frac{1}{4 \pi} \int_{-\infty}^{\infty} \frac{e^{i k x-\Omega^{(1)}_2(k) t}}{\Omega^{(1)}_2(k)} C^{(1)}_2(k,t) dk - \frac{1}{4 \pi} \int_{-\infty}^{\infty} \frac{e^{i k x-\Omega^{(1)}_1(k)t}}{\Omega^{(1)}_2(k)} C^{(1)}_1(k,t) dk.
\label{ec17_1}
\end{split}
\end{eqnarray}
On the other hand, proceeding similarly to \eqref{ec11_2} and applying the inverse Fourier transforms, we have
\begin{eqnarray}
\begin{split}
\Psi^{(2)}_1(x,t) = & \frac{1}{4 m_{2} \pi} \int_{-\infty}^{\infty}\left(\frac{\Omega^{(2)}_2(k) i+k}{\Omega^{(2)}_1(k)}\right) C^{(2)}_1(k,t) e^{i k x-\Omega^{(2)}_1(k) t} dk \\
& + \frac{1}{4 m_{2} \pi} \int_{-\infty}^{\infty}\left(\frac{\Omega^{(2)}_2(k) i-k}{\Omega^{(2)}_1(k)}\right) C^{(2)}_2(k,t) e^{i k x-\Omega^{(2)}_2(k)t} dk,
\label{ec16_2}
\end{split}
\end{eqnarray}
\begin{eqnarray}
\begin{split}
\Psi^{(2)}_2(x,t)=\frac{1}{4 \pi} \int_{-\infty}^{\infty} \frac{e^{i k x-\Omega^{(2)}_2(k) t}}{\Omega^{(2)}_2(k)} C^{(2)}_2(k,t) dk - \frac{1}{4 \pi} \int_{-\infty}^{\infty} \frac{e^{i k x-\Omega^{(2)}_1(k)t}}{\Omega^{(2)}_2(k)} C^{(2)}_1(k,t) dk.
\label{ec17_2}
\end{split}
\end{eqnarray}
\noindent The equations \eqref{ec16_1}-\eqref{ec17_2} are overdetermined because each contains boundary conditions, which are unnecessary. The first expression contains the boundary condition for $\Psi^{(1)}_2$, and the second expression contains the boundary condition for $\Psi^{(1)}_1$, analogously, we have the same for $\Psi^{(2)}_{1}$ and $\Psi^{(2)}_{2}$. All the aforementioned equations are superfluous to their respective expressions, and to remedy this, we utilize the discrete symmetry of the dispersion relations, $k \rightarrow -k$. Note, $\Omega^{(j)}_{1,2}(k)$ are invariant under $\nu(k)=-k$, since $\Omega^{(j)}_{1,2}(k) = \Omega^{(j)}_{1,2}(\nu(k))$, $j=1,2$. Applying this transformation to \eqref{ec16_1}, \eqref{ec17_1}, \eqref{ec16_2}, and \eqref{ec17_2}, we obtain 

\begin{eqnarray}
\begin{split}
\hat{\Psi}^{(1)}_1(-k,t) = & \left(\frac{\Omega^{(1)}_2(k)i-k}{2 m_{1} \Omega^{(1)}_1(k)}\right) C^{(1)}_1(-k,t) e^{-\Omega^{(1)}_1(k)t} \\
& + \left(\frac{\Omega^{(1)}_2(k)i+k}{2 m_{1} \Omega^{(1)}_1(k)}\right) C^{(1)}_2(-k,t) e^{-\Omega^{(1)}_2(k)t},
\label{ec18_1}
\end{split}
\end{eqnarray}

\begin{equation}
\hat{\Psi}^{(1)}_2(-k,t)=\frac{1}{2 \Omega^{(1)}_2(k)}\left(C^{(1)}_2(-k,t) e^{-\Omega^{(1)}_2(k)t}+C^{(1)}_1(-k,t) e^{-\Omega^{(1)}_1(k)t}\right),
\label{ec19_1}
\end{equation}

\begin{eqnarray}
\begin{split}
\hat{\Psi}^{(2)}_1(-k,t) = & \left(\frac{\Omega^{(2)}_2(k)i-k}{2 m_{2} \Omega^{(2)}_1(k)}\right) C^{(2)}_1(-k,t) e^{-\Omega^{(2)}_1(k)t} \\
& + \left(\frac{\Omega^{(2)}_2(k)i+k}{2 m_{2} \Omega^{(2)}_1(k)}\right) C^{(2)}_2(-k,t) e^{-\Omega^{(2)}_2(k)t},
\label{ec18_2}
\end{split}
\end{eqnarray}

\begin{equation}
\hat{\Psi}^{(2)}_2(-k,t)=\frac{1}{2 \Omega^{(2)}_2(k)}\left(C^{(2)}_2(-k,t) e^{-\Omega^{(2)}_2(k)t}+C^{(2)}_1(-k,t) e^{-\Omega^{(2)}_1(k)t}\right).
\label{ec19_2}
\end{equation}

\noindent Taking the inverse Fourier transform of \eqref{ec18_1}, \eqref{ec19_1}, \eqref{ec18_2}, and \eqref{ec19_2}, later, grouping all terms on one side yields expressions that sum to zero. Adding these expressions to \eqref{ec16_1}, \eqref{ec17_1}, \eqref{ec16_2}, and \eqref{ec17_2}, respectively. Finally, using the same argument used in Section $8$ of \cite{3} for each expression obtained and simplifying, we have the following result:

\begin{proposition}
The solution of the interface problem for the Dirac equation \eqref{ec1}-\eqref{ec5} for the massive case is given by  
\begin{eqnarray}
\begin{split}
& \Psi^{(1)}_{1}(x,t) = - \int_{-\infty}^{\infty} \frac{i k e^{ikx} (k \hat{\Psi}_{2,0}^{(1)}(-k) \sin(\alpha_{1}t) - m_{1} h_{0,1}^{(1)}(i\alpha_{1},t) (\sin(\alpha_{1}t)+i\cos(\alpha_{1}t)))}{2 \pi \alpha_{1}m_{1}} dk \\
& + \int_{-\infty}^{\infty}  \frac{e^{ikx} (-m_{1}\hat{\Psi}_{1,0}^{(1)}(-k)\cos(\alpha_{1}t)+m_{1}\hat{\Psi}_{1,0}^{(1)}(k)\cos(\alpha_{1}t)+i\alpha_{1}\hat{\Psi}_{2,0}^{(1)}(-k)\sin(\alpha_{1}t))}{2m_{1}\pi}dk \\
& - \frac{1}{2 m_{1} \pi} \int_{-\infty}^{\infty} \frac{i e^{ikx} \sin(\alpha_{1}t)(km_{1}\hat{\Psi}_{1,0}^{(1)}(-k) + k m_{1} \hat{\Psi}_{1,0}^{(1)}(k) + \alpha_{1}^{2}\hat{\Psi}_{2,0}^{(1)}(k))}{\alpha_{1}} dk \\
& + \frac{1}{2 m_{1} \pi} \int_{-\infty}^{\infty} \frac{ik e^{ikx} (k \hat{\Psi}_{2,0}^{(1)}(k) \sin(\alpha_{1}t) + h_{0,1}^{(1)}(-i\alpha_{1},t)(-m_{1} \sin(\alpha_{1}t) + im_{1} \cos(\alpha_{1}t)))}{\alpha_{1}} dk, 
\label{Psi_1_1_infinite}
\end{split}
\end{eqnarray}
and 
\begin{eqnarray}
\begin{split}
& \Psi^{(1)}_{2}(x,t) =  \int_{-\infty}^{\infty} \frac{e^{ikx} (k\hat{\Psi}_{2,0}^{(1)}(k)\cos(\alpha_{1}t) - \hat{\Psi}_{2,0}^{(1)}(-k)(k\cos(\alpha_{1}t) + i\alpha_{1}\sin(\alpha_{1}t)))}{2 \pi \alpha_{1}} dk \\
& + \int_{-\infty}^{\infty} \frac{e^{ikx} (m_{1}\hat{\Psi}_{1,0}^{(1)}(-k)\cos(\alpha_{1}t) -m_{1}\hat{\Psi}_{1,0}^{(1)}(k)\cos(\alpha_{1}t) + i\alpha_{1} \hat{\Psi}_{2,0}^{(1)}(k)\sin(\alpha_{1}t))}{2 \pi \alpha_{1}} dk \\
& + \frac{1}{2 \pi} \int_{-\infty}^{\infty} \frac{k e^{i(kx+\alpha_{1}t)} ( h_{0,1}^{(1)}(-i\alpha_{1},t) + h_{0,1}^{(1)}(i\alpha_{1},t))}{\alpha_{1}} dk,
\label{Psi_1_2_infinite}
\end{split}
\end{eqnarray}
for $x<0$, $t>0$, and $\alpha_{1}=\sqrt{k^2+m_{1}^{2}}$.\\
\begin{eqnarray}
\begin{split}
& \Psi^{(2)}_{1}(x,t) = - \int_{-\infty}^{\infty} \frac{i k e^{ikx} (k \hat{\Psi}_{2,0}^{(2)}(-k) \sin(\alpha_{2}t) + m_{2} h_{0,1}^{(2)}(i\alpha_{2},t) (\sin(\alpha_{2}t)+i\cos(\alpha_{2}t)))}{2 \pi \alpha_{2}m_{2}} dk \\
& + \int_{-\infty}^{\infty}  \frac{e^{ikx} (-m_{2}\hat{\Psi}_{1,0}^{(2)}(-k)\cos(\alpha_{2}t)+m_{2}\hat{\Psi}_{1,0}^{(2)}(k)\cos(\alpha_{2}t)+i\alpha_{2}\hat{\Psi}_{2,0}^{(2)}(-k)\sin(\alpha_{2}t))}{2m_{2}\pi}dk \\
& - \frac{1}{2 m_{2} \pi} \int_{-\infty}^{\infty} \frac{i e^{ikx} \sin(\alpha_{2}t)(km_{2}\hat{\Psi}_{1,0}^{(2)}(-k) + k m_{2} \hat{\Psi}_{1,0}^{(2)}(k) + \alpha_{2}^{2}\hat{\Psi}_{2,0}^{(2)}(k))}{\alpha_{2}} dk \\
& + \frac{1}{2 m_{2} \pi} \int_{-\infty}^{\infty} \frac{ik e^{ikx} (k \hat{\Psi}_{2,0}^{(2)}(k) \sin(\alpha_{2}t) - h_{0,1}^{(2)}(-i\alpha_{2},t)(-m_{2} \sin(\alpha_{2}t) + im_{2} \cos(\alpha_{2}t)))}{\alpha_{2}} dk, 
\label{Psi_2_1_infinite}
\end{split}
\end{eqnarray}
and 
\begin{eqnarray}
\begin{split}
& \Psi^{(2)}_{2}(x,t) =  \int_{-\infty}^{\infty} \frac{e^{ikx} (k\hat{\Psi}_{2,0}^{(2)}(k)\cos(\alpha_{2}t) - \hat{\Psi}_{2,0}^{(2)}(-k)(k\cos(\alpha_{2}t) + i\alpha_{2}\sin(\alpha_{2}t)))}{2 \pi \alpha_{2}} dk \\
& + \int_{-\infty}^{\infty} \frac{e^{ikx} (m_{2}\hat{\Psi}_{1,0}^{(2)}(-k)\cos(\alpha_{2}t) -m_{2}\hat{\Psi}_{1,0}^{(2)}(k)\cos(\alpha_{2}t) + i\alpha_{2} \hat{\Psi}_{2,0}^{(2)}(k)\sin(\alpha_{2}t))}{2 \pi \alpha_{2}} dk \\
& - \frac{1}{2 \pi} \int_{-\infty}^{\infty} \frac{k e^{i(kx+\alpha_{2}t)} ( h_{0,1}^{(2)}(-i\alpha_{2},t) + h_{0,1}^{(2)}(i\alpha_{2},t))}{\alpha_{2}} dk,
\label{Psi_2_2_infinite}
\end{split}
\end{eqnarray}
for $x>0$, $t>0$, and $\alpha_{2}=\sqrt{k^{2}+m_{2}^{2}}$.
\end{proposition}

\begin{remark}
It is important to note that \eqref{Psi_1_1_infinite}, \eqref{Psi_1_2_infinite}, \eqref{Psi_2_1_infinite}, and \eqref{Psi_2_2_infinite} are similar in form to the solutions of the Klein-Gordon equation in Section 7 of \cite{3}; this meshes well with the interpretation of the Dirac equation as a square root of a diagonal system of Klein-Gordon equations.
\end{remark}

\section{Solution of an interface problem for the Dirac equation on two finite domains}
\label{section_3} 

\quad We continue to generalize the solution for the interface problem for the Dirac equation by solving the system on two finite domains, $[-L, 0]$ and $[0; L]$, for $L>0$ finite. First, we solve the massless system and then the massive system. As before, the massive system has a more complex solution, but clear and expected consistency exists between the solutions of the two finite domains and the two semi-infinite domains. \\

\noindent The only change in the interface problem for the Dirac equation set-up from the previous section is the addition of another set of boundary conditions at $x=-L, \ L,$ and $x \in [-L, 0], [0, L]$, respectively. We seek four functions:
\begin{equation}
\begin{array}{ccc}
& \Psi_{1}^{(1)}(x,t), \Psi_{2}^{(1)}(x,t), \ x \in[-L, 0), t \in (0, T]; \\
& \Psi_{1}^{(2)}(x,t), \Psi_{2}^{(2)}(x,t), \ x \in(0, L], t \in (0, T], 
\label{Dirac_finite1}
\end{array}
\end{equation}
satisfying the equations
\begin{equation}
\begin{array}{ccc}
\left.
\begin{split}
i \partial_t \Psi_1^{(j)}(x,t)=-i \partial_x \Psi_1^{(j)}(x,t)+m_{j} \Psi_2^{(j)}(x,t),\ \textup{if }j=1, \ x \in [-L, 0); \\
i \partial_t \Psi_2^{(j)}(x,t)=i \partial_x \Psi_2^{(j)}(x,t)+m_{j} \Psi_1^{(j)}(x, t), \ \textup{if }j=2, \ x \in(0, L],
\label{Dirac_finite2}
\end{split} \right |
\ t \in(0, T],
\end{array}
\end{equation}
with initial conditions
\begin{equation}
\begin{array}{ccc}
\Psi_{1}^{(1)}(x,0)=\Psi_{1,0}^{(1)}(x), \ \Psi_{2}^{(1)}(x,0)=\Psi_{2,0}^{(1)}(x), \ x \in [-L,0); \\
\Psi_{1}^{(2)}(x,0)=\Psi_{1,0}^{(2)}(x), \ \Psi_{2}^{(2)}(x,0)=\Psi_{2,0}^{(2)}(x), \ x \in (0,L],
\label{Dirac_finite3}
\end{array}
\end{equation}
and interface continuity and boundary conditions
\begin{equation}
\begin{array}{ccc}
\Psi_{1}^{(1)}(0,t)=\Psi_{1}^{(2)}(0,t), \ \Psi_{2}^{(1)}(0,t)=\Psi_{2}^{(2)}(0,t), \ t \in (0,T]; \\
\Psi_{i}^{(1)}(-L,t)=\alpha_{i}^{(1)}(t), \ \Psi_{i}^{(2)}(L,t)=\alpha_{i}^{(2)}(t), \ t \in(0,T], \ i=1,2; 
\label{Dirac_finite4}
\end{array}
\end{equation}
where $T$ positive finite time value, $\alpha_{i}^{(1)}(t)$ are the boundary conditions for $\Psi_{i}^{(1)}(x,t)$ at $x= -L$, and $\alpha_{i}^{(2)}(t)$ are the boundary conditions for $\Psi_{i}^{(2)}(x,t)$ at $x=L$, for $i=1,2$.

\subsection{The massless system $(m_{j}=0, \ j=1,2)$ on two finite domains}
\label{section_3_1} 

All of the previous work in the last section for the massless and massive systems on two semi-infinite domains remains the same until the global relations are computed. The dispersion and local relations are not altered by the domain of the problem or the known boundary conditions; therefore, we begin the solution process for the massless Dirac equation on two finite domains at the computation of the global relations.

\noindent Analogously to Section \ref{section_2_1}, we have that the system given in \eqref{Dirac_finite1}-\eqref{Dirac_finite4} for the massless system is overdetermined. The UTM in the following steps illuminates that neither the left boundary condition for $\Psi^{(1)}_{1}$ nor the right boundary condition for $\Psi^{(2)}_{2}$ can be independently prescribed, and vice versa for $\Psi^{(2)}_{1}$ and $\Psi^{(2)}_{2}$, respectively. The initial conditions and the boundary conditions specify both of these conditions.

\noindent The same process is used to compute the global relations. We integrate the local relation for $\Psi_1^{(1)}(x,t)$ in the domain, $R_{1}=[-L,0] \times [0,T], \ L,T>0$ with $j=1$, we have
$$
\iint_{R_{1}}\left[e^{-i k x+\Omega_{1}(k) t} \Psi^{(1)}_1(x,t)\right]_t-\left[e^{-i k x+\Omega_{1}(k) t}\left(-\Psi^{(1)}_{1}(x,t)\right)\right]_x dtdx=0,
$$
and apply Green's Theorem to move the integration to the boundary, 
$$
\int_{\partial R_{1}}\left[e^{-i k x+\Omega_1(k)t} \Psi_1^{(1)}(x,t)\right] d x+\left[e^{-i k x+\Omega_{1}(k) t}\left(-\Psi^{(1)}_{1}(x,t)\right)\right] dt=0.
$$

\noindent Performing the integration on $R_{1}$ leads to the global relation, we have
\begin{equation}
\hat{\Psi}^{(1)}_1(k,0) + e^{i k L} B_{-L}\left(\Omega_{1}(k), t\right)-e^{\Omega_{1}(k)t} \hat{\Psi}^{(1)}_{1}(k,t)-B_{0}\left(\Omega_{1}(k),t\right)=0,
\label{ec29}
\end{equation}
where we define,
\begin{equation}
B_{-L}\left( \Omega_{j}(k),t \right)= \int_{0}^{t} e^{\Omega_{j}(k)s} \Psi^{(1)}_{j}(-L,s) ds = \int_{0}^{t} e^{\Omega_{j}(k)s} \alpha_{j}^{(1)} (t) ds, \ j=1,2;
\label{ec31}
\end{equation}

\begin{equation}
B_{0}\left(\Omega_{j}(k), t\right)=\int_0^t e^{\Omega_{j}(k)s} \Psi^{(1)}_j(0,s) ds=\int_0^t e^{\Omega_{j}(k)s} \Psi^{(2)}_j(0,s)ds, \ j=1,2.
\label{ec30}
\end{equation}
\noindent Following the same procedure for the local relation of $\Psi^{(1)}_{2}(x,t)$ gives the second global relation,
\begin{equation}
\hat{\Psi}^{(1)}_2(k,0) - e^{ikL} B_{-L}\left(\Omega_{2}(k), t\right)-e^{\Omega_{2}(k) t} \hat{\Psi}^{(1)}_{2}(k,t) + B_0\left(\Omega_{2}(k),t\right)=0.
\label{ec32}
\end{equation}

\noindent Now, we solve for $\hat{\Psi}^{(1)}_1(k, t)$ and $\hat{\Psi}^{(1)}_2(k,t)$ in \eqref{ec29} and \eqref{ec32} and apply the inverse Fourier transform. Solving for the Fourier transforms, we obtain
$$
\begin{aligned}
& \hat{\Psi}^{(1)}_1(k,t)=e^{-\Omega_{1}(k)t}\left(\hat{\Psi}^{(1)}_{1}(k,0) + e^{ikL} B_{-L}\left(\Omega_{1}(k),t\right)-B_{0}\left(\Omega_{1}(k),t\right)\right), \\
& \hat{\Psi}^{(1)}_2(k,t)=e^{-\Omega_{2}(k)t}\left(\hat{\Psi}^{(1)}_{2}(k,0) - e^{-ikt} B_{-L}\left(\Omega_{2}(k),t\right)+B_{0}\left(\Omega_{2}(k),t\right)\right),
\end{aligned}
$$
and taking the inverse Fourier transform of these expressions leads to the following
\begin{equation}
\Psi^{(1)}_{1}(x,t)=\frac{1}{2 \pi} \int_{-\infty}^{\infty} e^{i k(x-t)}\left(\hat{\Psi}^{(1)}_{1,0}(k)+e^{i k L} B_{-L}\left(\Omega_{1}(k),t\right)-B_{0}\left(\Omega_{1}(k),t\right)\right) dk, 
\label{ec33}
\end{equation}
\begin{equation}
\Psi^{(1)}_{2}(x,t)=\frac{1}{2 \pi} \int_{-\infty}^{\infty} e^{ik(x+t)}\left(\hat{\Psi}^{(1)}_{2,0}(k)-e^{i k L} B_{-L}\left(\Omega_{2}(k),t\right)+B_{0}\left(\Omega_{2}(k),t\right)\right)dk.
\label{ec34}
\end{equation}

\noindent For the other part, integrating the local relation \eqref{ecuacion6} on the domain $R_{2} = [0,L] \times [0,T], \ L, T > 0$ when $j=2$, and using Green's Theorem,  we obtain
\begin{eqnarray}
\begin{split}
& \hat{\Psi}^{(2)}_{1}(k,0)-e^{-ikL} B_{L}\left(\Omega_{1}(k), t\right)-e^{\Omega_{1}(k) t} \hat{\Psi}^{(2)}_{1}(k,t)+ B_{0}\left(\Omega_{1}(k), t\right)=0, \\
& \hat{\Psi}^{(2)}_{2}(k,0)+e^{-ikL} B_{L}\left(\Omega_{2}(k), t\right)-e^{\Omega_{2}(k) t} \hat{\Psi}^{(2)}_{2}(k,t)-B_{0}\left(\Omega_{2}(k), t\right)=0.
\label{ecuacion8_finito}
\end{split}
\end{eqnarray}
where we define
$$
B_{L}\left(\Omega_{j}(k), t\right) = \int_{0}^{t} e^{\Omega_{j}(k)s} \Psi^{(2)}_{j}(L,s)ds = \int_{0}^{t} e^{\Omega_{j}(k)s} \alpha_{j}^{(2)}(t) ds, \ j=1,2.
$$

\noindent The next step is to utilize the inverse Fourier transform to obtain expressions for $\Psi_1^{(2)}(x,t)$ and $\Psi_2^{(2)}(x, t)$ from \eqref{ecuacion8_finito}. Multiplying the first equation in \eqref{ecuacion8_finito} by $e^{-\Omega_{1}(k)t}$ and the second by $e^{-\Omega_{2}(k)t}$, then algebraically solving for the Fourier transforms yields
$$
\begin{aligned}
& \hat{\Psi}^{(2)}_{1}(k,t)=e^{-\Omega_{1}(k) t}\left(\hat{\Psi}^{(2)}_{1}(k,0)-e^{-i k L} B_{L}\left(\Omega_{1}(k),t\right)+B_{0}\left(\Omega_{1}(k), t\right)\right), \\
& \hat{\Psi}^{(2)}_{2}(k,t)=e^{-\Omega_{2}(k) t}\left(\hat{\Psi}^{(2)}_{2}(k,0)+e^{-i k L} B_{L}\left(\Omega_{2}(k),t\right)-B_{0}\left(\Omega_{2}(k), t\right)\right).
\end{aligned}
$$

\noindent Taking the inverse Fourier transforms of the above expressions and substituting them into the dispersion relations \eqref{ecuacion3}, we have the following
\begin{equation}
\Psi^{(2)}_{1}(x,t)=\frac{1}{2 \pi} \int_{-\infty}^{\infty} e^{i k(x-t)}\left(\hat{\Psi}^{(2)}_{1,0}(k)-e^{-i k L} B_{L}\left(\Omega_{1}(k),t\right)+B_{0}\left(\Omega_{1}(k), t\right)\right) dk, 
\label{ecuacion9_1_finito}
\end{equation}
\begin{equation}
\Psi^{(2)}_{2}(x,t)=\frac{1}{2 \pi} \int_{-\infty}^{\infty} e^{i k(x+t)}\left(\hat{\Psi}^{(2)}_{2,0}(k)+e^{-i k L} B_{L}\left(\Omega_{2}(k),t\right)-B_{0}\left(\Omega_{2}(k), t\right)\right) dk.
\label{ecuacion9_2_finito}
\end{equation}

\noindent As in the two semi-infinite domains in the massless case, we can simplify \eqref{ec33} and \eqref{ec34} about the boundary conditions to remove the unnecessary right boundary condition of $\Psi^{(1)}_{1}$ and the left boundary condition of $\Psi^{(1)}_{2}$. Expanding \eqref{ec33} and \eqref{ec34} with definitions \eqref{ec31} and \eqref{ec30},
\begin{eqnarray}
\begin{split}
\Psi^{(1)}_{1}(x,t) = & \frac{1}{2 \pi} \int_{-\infty}^{\infty} e^{i k(x-t)} \hat{\Psi}^{(1)}_{1,0}(k) dk + \frac{1}{2 \pi} \int_{-\infty}^{\infty} \int_0^t e^{i k(x-t+L+s)} \alpha^{(1)}_{1}(s) ds dk \\
& - \frac{1}{2 \pi} \int_{-\infty}^{\infty} \int_0^t e^{i k(x-t+s)} \Psi_{1}^{(1)}(0,s) dsdk,
\label{ec35}
\end{split}
\end{eqnarray}
\begin{eqnarray}
\begin{split} 
\Psi^{(1)}_{2}(x,t) = & \frac{1}{2 \pi} \int_{-\infty}^{\infty} e^{i k(x+t)} \hat{\Psi}^{(1)}_{2,0}(k) dk - \frac{1}{2 \pi} \int_{-\infty}^{\infty} \int_0^t e^{i k(x+t+L-s)} \alpha^{(1)}_{2}(s) dsdk \\
& + \frac{1}{2 \pi} \int_{-\infty}^{\infty} \int_0^t e^{i k(x+t-s)} \Psi^{(1)}_{2}(0,s) dsdk.
\label{ec36}
\end{split}
\end{eqnarray}
\noindent Proceeding similarly with \eqref{ecuacion9_1_finito} and \eqref{ecuacion9_2_finito}, we found $\Psi_1^{(2)}(x,t)$ and $\Psi_2^{(2)}(x,t)$, respectively. Thus, we have
\begin{eqnarray}
\begin{split}
\Psi^{(2)}_1(x,t) = & \frac{1}{2 \pi} \int_{-\infty}^{\infty} e^{i k(x-t)} \hat{\Psi}^{(2)}_{1,0}(k) dk-\frac{1}{2 \pi} \int_{-\infty}^{\infty} \int_0^t e^{i k(x-t-L+s)} \alpha_{1}^{(2)}(s) dsdk \\
& + \frac{1}{2 \pi} \int_{-\infty}^{\infty} \int_0^t e^{i k(x-t+s)} \Psi^{(2)}_1(0,s) dsdk,
\label{ec37}
\end{split}
\end{eqnarray}

\begin{eqnarray}
\begin{split}
\Psi^{(2)}_2(x,t) = & \frac{1}{2 \pi} \int_{-\infty}^{\infty} e^{i k(x+t)} \hat{\Psi}^{(2)}_{2,0}(k) dk + \frac{1}{2 \pi} \int_{-\infty}^{\infty} \int_0^t e^{i k(x+t-L-s)} \alpha^{(2)}_2(s) dsdk \\
& - \frac{1}{2 \pi} \int_{-\infty}^{\infty} \int_0^t e^{i k(x+t-s)} \Psi^{(2)}_2(0,s) dsdk.
\label{ec64}
\end{split}
\end{eqnarray}

\noindent Simplifying \eqref{ec35} first, since $x>-L$ and $t>s$, this implies $x-t+L+s>0$. Using a contour similar to $\Gamma_{C}$ shown in Figure \ref{Figura_2}, except with $\textup{Arc}_{C}$ in the upper half plane in $\mathbb{C}$, we can apply Cauchy's Theorem and Jordan's Lemma to show the integral vanishes due to the exponential decay and analytic properties of the integrand. The new contour is shown below in Figure \ref{Figura_3}.

\begin{figure}[h]
\centering
\includegraphics[height=5cm,width=9cm]{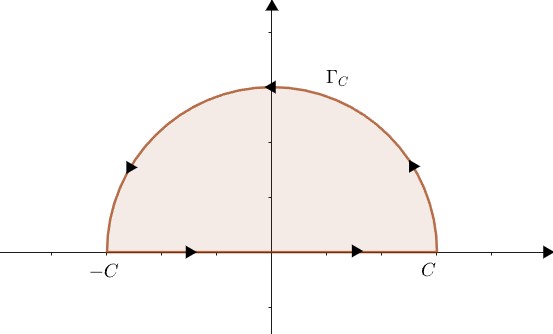}
\caption{Contour utilized to simplify \eqref{ec35}, which is similar to contour $\Gamma_{C}$ in the previous section.}
\label{Figura_3}
\end{figure}

\noindent The last integral in \eqref{ec35} is a previously treated case in the Subsection \ref{section_2_1}, so applying the inverse Fourier transform,
$$
\Psi^{(1)}_1(x,t)= \begin{cases}\Psi^{(1)}_{1,0}(x-t), & x - t \geq -L; \\ \Psi^{(1)}_{1}(-L,x-t+L), & x - t < -L.\end{cases}
$$

\noindent In \eqref{ec36}, the last integral was shown to vanish in Subsection \ref{section_2_1} by applying the contour $\Gamma_{C}$. The second integral in \eqref{ec36} is not signed because $x-s<0$ and $L+t>0$, so we use the definition of the delta distribution to obtain
$$
- \frac{1}{2 \pi} \int_{-\infty}^{\infty} \int_0^t e^{i k(x+t+L-s)} \alpha^{(1)}_2(s) dsdk = - \alpha^{(1)}_2(x+t+L).
$$

\noindent When $x+t+L>0$, we have $\alpha^{(1)}_{2}(x+t+L)=0$ because the boundary condition is only defined for positive time values. Therefore, applying the inverse Fourier transform to these cases,
$$
\Psi^{(1)}_2(x,t)= \begin{cases}\Psi^{(1)}_{2,0}(x+t), & -x \leq t; \\ \Psi_{2}^{(1)}(0,x+t), & -x > t,\end{cases}
$$
where $\Psi^{(1)}_{2,0}(x+t)=0$ when $x+t>-L$ because these arguments are outside the domain where initial conditions are specified. \\

\noindent Thus, we have the solution for the interface problem for the massless Dirac equation on the finite interval $[-L,0]$. On the other hand, proceeding analogously to how $\Psi^{(1)}_{1}(x,t)$ and $\Psi^{(1)}_{2}(x,t)$ were found and using similar arguments, we obtain the solution for the interface problem for the massless Dirac equation on $[0,L]$, that is, we find $\Psi^{(2)}_{1}(x,t)$ and $\Psi^{(2)}_{2}(x,t)$. Therefore, we have the following result:

\begin{proposition}
The solution of the interface problem for the Dirac equation \eqref{Dirac_finite1}-\eqref{Dirac_finite4} for the massless case is given by  
\begin{eqnarray}
\begin{split}
& \Psi^{(1)}_1(x,t)= \begin{cases}\Psi_{1,0}^{(1)}(x-t), & x - t \geq -L; \\
\Psi_{1}^{(1)}(-L,x - t + L), & x - t < -L,\end{cases} 
\label{ec39_a}
\end{split}
\end{eqnarray}
\noindent and
\begin{eqnarray}
\begin{split}
& \Psi_{2}^{(1)}(x,t)= \begin{cases}\Psi_{2,0}^{(1)}(x+t), & -x \leq t; \\
\Psi_{2}^{(1)}(0, x+t), & -x > t, \end{cases}
\label{ec39_b}
\end{split}
\end{eqnarray}
for $x<0$,\ $t>0$. \\
\begin{eqnarray}
\begin{split}
& \Psi^{(2)}_1(x, t)= \begin{cases}\Psi^{(2)}_{1,0}(x-t), & x \geq t; \\
\Psi^{(2)}_1(0,t-x), & x<t, \end{cases} 
\label{ec40_a}
\end{split}
\end{eqnarray}
\noindent and
\begin{eqnarray}
\begin{split}
& \Psi_{2}^{(2)}(x,t)= \begin{cases}\Psi_{2,0}^{(2)}(x+t), & x+t \leq L; \\
\Psi_{2}^{(2)}(L, x+t-L), & x+t>L,\end{cases}
\label{ec40_b}
\end{split}
\end{eqnarray}
for $x>0$,\ $t>0$.
\end{proposition}

\subsection{The massive system $(m_{j}>0, \ j=1,2)$ on two finite domains}
\label{section_3_2}

For the massive system, we again start with the global relations because the dispersion and local relations are not affected by restricting the domain of the PDEs to a finite region in the plane. These properties are only dictated by local information, and they are preserved by the global information in the PDEs used.

\noindent Analogous to the massive case on two semi-infinite domains, the system as stated for the massive case on two finite domains is not overdetermined. However, the final solutions computed via the UTM will show that the solution for $\Psi^{(1)}_{1}$ does not depend on the right boundary condition of $\Psi^{(1)}_{2}$, and the solution of $\Psi^{(1)}_{2}$ does not depend on the right boundary condition of $\Psi^{(1)}_{1}$. Instead for $\Psi^{(2)}_{1}$ and $\Psi^{(2)}_{2}$, we have that the solution for $\Psi^{(2)}_{1}$ does not depend on the left boundary condition of $\Psi^{(2)}_{2}$, and the solution of $\Psi^{(2)}_{2}$ does not depend on the left boundary condition of $\Psi^{(2)}_{1}$. \\

\noindent Using the same process to compute the global relation as before, we integrate the local relations over the domain, $R_{1}=[-L, 0] \times [0, T]$ and $R_{2}=[0, L] \times [0, T]$, $L,T>0$, respectively. Thus, we have
\begin{eqnarray}
\begin{split}
& \iint_{R_{1}}\left(e^{-i k x+\Omega^{(1)}_j(k) t}\left[(i m_{1}) \Psi^{(1)}_1(x,t)+\left(\Omega^{(1)}_j(k)-i k\right) \Psi^{(1)}_2(x,t)\right]\right)_t \\
& - \left(e^{-i k x+\Omega^{(1)}_{j}(k)t}\left[(-im_{1}) \Psi^{(1)}_{1}(x,t)+\left(\Omega^{(1)}_j(k)-i k\right) \Psi^{(1)}_{2}(x,t)\right]\right)_x dxdt=0,\ j=1,2;
\label{RL1_section3}
\end{split}
\end{eqnarray}
\begin{eqnarray}
\begin{split}
& \iint_{R_{2}}\left(e^{-i k x+\Omega^{(2)}_j(k) t}\left[(i m_{2}) \Psi^{(2)}_1(x,t)+\left(\Omega^{(2)}_j(k)-i k\right) \Psi^{(2)}_2(x,t)\right]\right)_t \\
& - \left(e^{-i k x+\Omega^{(2)}_{j}(k)t}\left[(-im_{2}) \Psi^{(2)}_{1}(x,t)+\left(\Omega^{(2)}_j(k)-i k\right) \Psi^{(2)}_{2}(x,t)\right]\right)_x dxdt=0,\ j=1,2;
\label{RL2_section3}
\end{split}
\end{eqnarray}
and applying Green's Theorem to move the integration to the boundary in \eqref{RL1_section3} and \eqref{RL2_section3}, respectively, we obtain
\begin{eqnarray}
\begin{split}
& \int_{\partial R_{1}}\left(e^{-i k x \Omega^{(1)}_{j}(k) t}\left[(i m_{1}) \Psi^{(1)}_{1}(x,t)+\left(\Omega^{(1)}_{j}(k)-i k\right) \Psi^{(1)}_2(x,t)\right]\right) dx  \\
& + \left(e^{-i k x+\Omega^{(1)}_{j}(k)t}\left[(-im_{1}) \Psi^{(1)}_{1}(x,t)+\left(\Omega^{(1)}_{j}(k)-i k\right) \Psi^{(1)}_{2}(x,t)\right]\right)dt=0, \ j=1,2;
\label{RG1_section3}
\end{split}
\end{eqnarray}
\begin{eqnarray}
\begin{split}
& \int_{\partial R_{2}}\left(e^{-i k x \Omega^{(2)}_j(k) t}\left[(i m_{2}) \Psi^{(2)}_1(x,t)+\left(\Omega^{(2)}_j(k)-i k\right) \Psi^{(2)}_{2}(x,t)\right]\right) dx  \\
& + \left(e^{-i k x+\Omega^{(2)}_{j}(k)t}\left[(-im_{2}) \Psi^{(2)}_{1}(x,t)+\left(\Omega^{(2)}_{j}(k)-i k\right) \Psi^{(2)}_{2}(x,t)\right]\right)dt=0, \ j=1,2.
\label{RG2_section3}
\end{split}
\end{eqnarray}

\noindent Performing the necessary line integrals gives the following global relations for \eqref{RG1_section3} and \eqref{RG2_section3}, respectively. We have
\begin{eqnarray}
\begin{split}
& (im_{1}) \hat{\Psi}^{(1)}_1(k,0) + \left(\Omega^{(1)}_j(k)-i k\right) \hat{\Psi}^{(1)}_2(k,0) \\
& +(im_{1}) e^{ikL} B^{(1)}_{-L,1}\left(\Omega^{(1)}_j(k),t\right)-\left(\Omega^{(1)}_j(k)-i k\right) e^{ikL} B^{(1)}_{-L,2}\left(\Omega^{(1)}_j(k),t\right) \\
& +(-im_{1}) e^{\Omega^{(1)}_{j}(k) t} \hat{\Psi}^{(1)}_1(k,t)+\left(i k-\Omega^{(1)}_j(k)\right) e^{\Omega^{(1)}_j(k)t} \hat{\Psi}^{(1)}_2(k,t) \\
& +(-im_{1}) B^{(1)}_{0,1}\left(\Omega^{(1)}_j(k),t\right) - \left(i k - \Omega^{(1)}_j(k)\right) B^{(1)}_{0,2}\left(\Omega^{(1)}_j(k), t\right)=0, \ j=1,2,
\label{ec38}
\end{split}
\end{eqnarray}
where we define,
\begin{equation}
B^{(1)}_{-L,i}\left(\Omega^{(1)}_j(k),t\right)=\int_0^t e^{\Omega^{(1)}_j(k)s} \Psi^{(1)}_i(-L,s) ds=\int_0^t e^{\Omega^{(1)}_j(k)s} \alpha^{(1)}_i(s) ds, \ i,j=1,2.
\label{ec39}
\end{equation}
\begin{equation}
B^{(1)}_{0,i}\left(\Omega^{(1)}_j(k),t\right)=\int_0^t e^{\Omega^{(1)}_j(k)s} \Psi^{(1)}_i(0,s) ds=\int_0^t e^{\Omega^{(1)}_j(k)s} \Psi^{(2)}_i(0,s) ds, \ i,j=1,2.
\label{ec40}
\end{equation}

\noindent and
\begin{eqnarray}
\begin{split}
& (im_{2}) \hat{\Psi}^{(2)}_1(k,0) + \left(\Omega^{(2)}_{j}(k)-i k\right) \hat{\Psi}^{(2)}_{2}(k,0) \\
& +(-im_{2}) e^{-ikL} B^{(2)}_{L,1}\left(\Omega^{(2)}_{j}(k),t\right)+\left(\Omega^{(2)}_{j}(k)-i k\right) e^{-ikL} B^{(2)}_{L,2}\left(\Omega^{(2)}_{j}(k),t\right) \\
& +(-im_{2}) e^{\Omega^{(2)}_{j}(k) t} \hat{\Psi}^{(2)}_{1}(k,t)+\left(i k-\Omega^{(2)}_{j}(k)\right) e^{\Omega^{(2)}_{j}(k)t} \hat{\Psi}^{(2)}_{2}(k,t) \\
& +(im_{2}) B^{(2)}_{0,1}\left(\Omega^{(2)}_{j}(k),t\right)+\left(i k - \Omega^{(2)}_{j}(k)\right) B^{(2)}_{0,2}\left(\Omega^{(2)}_{j}(k), t\right)=0, \ j=1,2,
\label{ec41}
\end{split}
\end{eqnarray}
where we define
\begin{equation}
B^{(2)}_{0,i}\left(\Omega^{(2)}_{j}(k),t\right)=\int_{0}^{t} e^{\Omega^{(2)}_{j}(k)s} \Psi^{(2)}_i(0,s) ds=\int_{0}^{t} e^{\Omega^{(2)}_{j}(k)s} \Psi^{(1)}_{i}(0,s) ds, \ i,j=1,2;
\label{ec42}
\end{equation}
\begin{equation}
B^{(2)}_{L,i}\left(\Omega^{(2)}_{j}(k),t\right)=\int_{0}^{t} e^{\Omega^{(2)}_{j}(k)s} \Psi^{(2)}_{i}(L,s) ds=\int_{0}^{t} e^{\Omega^{(2)}_{j}(k)s} \alpha^{(2)}_{i}(s) ds, \ i,j=1,2.
\label{ec43}
\end{equation}

\noindent Comparing these global relations to the ones for the massive case on two semi-infinite domains, many similarities are noted with the addition of two extra terms and the new definitions for the integrals on the left and right boundaries. These similarities ease the algebra necessary to obtain valid expressions for $\Psi^{(1)}_{1}(x,t)$, $\Psi^{(1)}_{2}(x,t)$, $\Psi^{(2)}_{1}(x,t)$, and $\Psi^{(2)}_{2}(x,t)$, respectively. \\

\noindent Now, we use the global relations and solve for $\hat{\Psi}^{(1)}_{1}(k,t)$, $\hat{\Psi}^{(1)}_{2}(k,t)$, $\hat{\Psi}^{(2)}_{1}(k,t)$, and $\Psi^{(2)}_{2}(x,t)$, to apply the inverse Fourier transform. This is achieved on two semi-infinite domains by multiplying the global relations by different exponentials and adding and subtracting the respective expressions. Due to the symmetry between the two semi-infinite domains and two finite domains cases, the same steps as before can be facts, but new definitions for $C^{(1)}_{j}(k,t)$ in \eqref{ec33} and \eqref{ec34}, with $j=1,2$; and for $C^{(2)}_{j}(k,t)$ in \eqref{ecuacion9_1_finito} and \eqref{ecuacion9_2_finito}, with $j=1,2$, it is necessary to account for the additional terms left from integrating over the two finite domains, respectively. For two finite domains, the new definitions for these expressions are the following: 
\begin{eqnarray*}
\begin{split}
C^{(1)}_{j}(k,t)= & (i m_{1}) \hat{\Psi}^{(1)}_{1,0}(k)+\left(\Omega^{(1)}_j(k)-i k\right) \hat{\Psi}^{(1)}_{2,0}(k) \\
& +(im_{1}) e^{ikL} B^{(1)}_{-L,1}\left(\Omega^{(1)}_j(k),t\right) - \left(\Omega^{(1)}_j(k)-i k\right) e^{ikL} B^{(1)}_{-L,2}\left(\Omega^{(1)}_{j}(k),t\right) \\
& +(-im_{1}) B^{(1)}_{0,1}\left(\Omega^{(1)}_{j}(k),t\right) - \left(ik-\Omega^{(1)}_{j}(k)\right) B^{(1)}_{0,2}\left(\Omega^{(1)}_{j}(k),t\right), \ j = 1,2;
\end{split}
\end{eqnarray*}
\begin{eqnarray*}
\begin{split}
C^{(2)}_{j}(k,t)= & (i m_{2}) \hat{\Psi}^{(2)}_{1,0}(k)+\left(\Omega^{(2)}_{j}(k)-i k\right) \hat{\Psi}^{(2)}_{2,0}(k) \\
& -(im_{2}) e^{-ikL} B^{(2)}_{L,1}\left(\Omega^{(2)}_{j}(k),t\right)+\left(\Omega^{(2)}_{j}(k)-ik\right) e^{-ikL} B^{(2)}_{L,2}\left(\Omega^{(2)}_{j}(k),t\right) \\
& +(im_{2}) B^{(2)}_{0,1}\left(\Omega^{(2)}_{j}(k),t\right)+\left(ik-\Omega^{(2)}_{j}(k)\right) B^{(2)}_{0,2}\left(\Omega^{(2)}_{j}(k),t\right), \ j = 1,2.
\end{split}
\end{eqnarray*}

\noindent Analogously to two semi-infinite domains, applying discrete symmetry $k \rightarrow-k$ enables us to eliminate unnecessary boundary terms and obtain the solution on the finite interval $[-L,0]$. Working through the same argument, we obtain our solutions $\Psi^{(1)}_{1}(x,t)$ and $\Psi^{(1)}_{2}(x,t)$. We were able to remove the right boundary condition of $\Psi^{(1)}_{2}$, contained in $B^{(1)}_{0,2}$ from the expression obtained for $\Psi^{(1)}_{1}(x,t)$ and the right boundary condition of $\Psi^{(1)}_{1}$, contained in $B^{(1)}_{0,1}$ from the expression obtained for $\Psi^{(1)}_{2}(x,t)$. Thus, we obtain the solution to the interface problem for the Dirac equation for the massive case on $[-L,0]$. For the other part, using arguments similar to those used, we obtain expressions for $\Psi_{1}^{(2)}(x,t)$ and $\Psi_{2}^{(2)}(x,t)$, then applying the arguments given in Subsection \ref{section_2_2}, we have the solution to the interface problem for the Dirac equation for the massive case on $[0,L]$. Therefore, we obtain the following result:

\begin{proposition}
The solution of the interface problem for the Dirac equation \eqref{Dirac_finite1}-\eqref{Dirac_finite4} for the massive case is given by 
\begin{eqnarray}
\begin{split}
& \Psi_{1}^{(1)}(x,t) = \frac{1}{4 m_{1} \pi} \int_{-\infty}^{\infty} \frac{e^{ikx}\cos(\alpha_{1}t)m_{1}\alpha_{1}i(2\hat{\Psi}_{1,0}^{(1)}(k)-2\hat{\Psi}_{1,0}^{(1)}(-k) - B_{-L,1}^{(1)}(i\alpha_{1},t)e^{-ikL})}{\alpha_{1}i} dk \\
& + \int_{-\infty}^{\infty} \frac{e^{ikx}\cos(\alpha_{1}t)m_{1}\alpha_{1}i( B_{-L,1}^{(1)}(i\alpha_{1},t)e^{ikL} - B_{-L,1}^{(1)}(-i\alpha_{1},t)(e^{-ikL} - e^{ikL}) )}{4 \alpha_{1} m_{1} \pi i}dk \\
& - \int_{-\infty}^{\infty} \frac{e^{ikx}m_{1}\alpha_{1}\sin(\alpha_{1}t)(B_{-L,1}^{(1)}(i\alpha_{1},t)(e^{-ikL} - e^{ikL}) + B_{-L,1}^{(1)}(-i\alpha_{1},t) (e^{ikL} - e^{-ikL}) )}{4 \alpha_{1} m_{1} \pi i} dk \\
%%%%
& + \int_{-\infty}^{\infty} \frac{e^{ikx}m_{1}k\cos(\alpha_{1}t)i(B_{-L,1}^{(2)}(i\alpha_{1},t)(e^{-ikL} + e^{ikL}) - B_{-L,1}^{(1)}(-i\alpha_{1},t)(e^{-ikL} + e^{ikL}))}{4 \alpha_{1} m_{1} \pi i} dk \\
& + \int_{-\infty}^{\infty} \frac{2e^{ikx}\sin(\alpha_{2}t)\alpha_{2}^{2}(\hat{\Psi}_{2,0}^{(2)}(k) - \hat{\Psi}_{2,0}^{(2)}(-k)) + 2e^{ikx}\sin(\alpha_{2}t)k^{2}(\hat{\Psi}_{2,0}^{(2)}(-k) - \hat{\Psi}_{2,0}^{(2)}(k)) }{4 \alpha_{2} m_{2} \pi i} dk \\
& + \int_{-\infty}^{\infty} \frac{e^{ikx}\cos(\alpha_{1}t)(\alpha_{1}^{2}-k^{2})i( B_{-L,2}^{(1)}(-i\alpha_{1},t) - B_{-L,2}^{(1)}(i\alpha_{1},t) )(e^{ikL} - e^{-ikL})}{4 \alpha_{1} m_{1} \pi i}dk \\
& - \int_{-\infty}^{\infty} \frac{e^{ikx}\sin(\alpha_{1}t)(\alpha_{1}^{2} - k^{2})(B_{-L,2}^{(1)}(i\alpha_{1},t) + B_{-L,2}^{(1)}(-i\alpha_{1},t) )(e^{ikL} - e^{-ikL})}{4 \alpha_{1} m_{1} \pi i}dk \\
& + \int_{-\infty}^{\infty} \frac{2e^{ikx}\cos(\alpha_{1}t)m_{1}ki( B_{0,1}^{(1)}(-i\alpha_{1},t) - B_{0,1}^{(1)}(i\alpha_{1},t) )}{4 \alpha_{1} m_{1} \pi i}dk \\
& - \int_{-\infty}^{\infty} \frac{2e^{ikx}\sin(\alpha_{1}t)m_{1}k(B_{0,1}^{(1)}(i\alpha_{1},t) + B_{0,1}^{(1)}(-i\alpha_{1},t) )}{4 \alpha_{1} m_{1} \pi i}dk \\
& + \int_{-\infty}^{\infty} \frac{e^{ikx}m_{1}k\sin(\alpha_{1}t)(B_{-L,1}^{(1)}(i\alpha_{1},t) + B_{-L,1}^{(1)}(-i\alpha_{1},t)) (e^{-ikL} + e^{ikL})}{4 \alpha_{1} m_{1} \pi i} dk \\
& + \int_{-\infty}^{\infty} \frac{e^{ikx}m_{1}k\sin(\alpha_{1}t)(2\hat{\Psi}_{1,0}^{(1)}(k) + 2\hat{\Psi}_{1,0}^{(1)}(-k))}{4 \alpha_{1} m_{1} \pi i} dk,
\label{Psi_1_1}
\end{split}
\end{eqnarray}
and 
\begin{eqnarray}
\begin{split}
& \Psi_{2}^{(1)}(x,t) =  \int_{-\infty}^{\infty} \frac{e^{ikx}\cos(\alpha_{1}t)m_{1}i(\hat{\Psi}_{1,0}^{(1)}(k) - \hat{\Psi}_{1,0}^{(1)}(-k)) }{2 \pi \alpha_{1} i} dk \\
& + \int_{-\infty}^{\infty} \frac{ e^{ikx}\cos(\alpha_{1}t)\alpha_{1}i(\hat{\Psi}_{2,0}^{(1)}(-k) - \hat{\Psi}_{2,0}^{(1)}(k))}{2 \pi \alpha_{1} i} dk \\
& - \int_{-\infty}^{\infty} \frac{e^{ikx} \cos(\alpha_{1}t)ki \hat{\Psi}_{2,0}^{(1)}(k) + e^{ikx} \sin(\alpha_{1}t)k \hat{\Psi}_{2,0}^{(1)}(-k)}{2 \pi \alpha_{1} i} dk \\
& - \int_{-\infty}^{\infty} \frac{e^{ikx} \cos(\alpha_{1}t) (m_{1} B_{-L,1}^{(1)}(-i\alpha_{1},t) + \alpha_{1} B_{-L,2}^{(1)}(-i\alpha_{1},t))i(e^{-ikL} - e^{ikL})}{2 \pi \alpha_{1} i} dk \\
& + \int_{-\infty}^{\infty} \frac{e^{ikx} \cos(\alpha_{1}t) k i B_{-L,2}^{(1)}(-i\alpha_{1},t)(e^{-ikL} + e^{ikL} - 2)}{2 \pi \alpha_{1} i} dk,
\label{Psi_1_2}
\end{split}
\end{eqnarray}
for $x<0$, $t>0$, and $\alpha_{1}=\sqrt{k^{2}+m_{1}^{2}}$.
\begin{eqnarray}
\begin{split}
& \Psi_{1}^{(2)}(x,t) = \frac{1}{4 m_{2} \pi} \int_{-\infty}^{\infty} \frac{e^{ikx}\cos(\alpha_{2}t)m_{2}\alpha_{2}i(2\hat{\Psi}_{1,0}^{(2)}(k)-2\hat{\Psi}_{1,0}^{(2)}(-k)+B_{L,1}^{(2)}(i\alpha_{2},t)e^{ikL})}{\alpha_{2}i} dk \\
& + \int_{-\infty}^{\infty} \frac{e^{ikx}\cos(\alpha_{2}t)m_{2}\alpha_{2}i(B_{L,1}^{(2)}(-i\alpha_{2},t)(e^{ikL} - e^{-ikL})- B_{L,1}^{(2)}(i\alpha_{2},t)e^{-ikL} )}{4 \alpha_{2} m_{2} \pi i}dk \\
& + \int_{-\infty}^{\infty} \frac{e^{ikx}m_{2}\alpha_{2}\sin(\alpha_{2}t)(B_{L,1}^{(2)}(i\alpha_{2},t)(e^{ikL} - e^{-ikL}) + B_{L,1}^{(2)}(-i\alpha_{2},t) (e^{-ikL} - e^{ikL}) )}{4 \alpha_{2} m_{2} \pi i} dk \\
%%%%
& + \int_{-\infty}^{\infty} \frac{e^{ikx}m_{2}k\cos(\alpha_{2}t)i(B_{L,1}^{(2)}(-i\alpha_{2},t)(e^{ikL} + e^{-ikL}) - B_{L,1}^{(2)}(i\alpha_{2},t)(e^{ikL} + e^{-ikL}))}{4 \alpha_{2} m_{2} \pi i} dk \\
& + \int_{-\infty}^{\infty} \frac{2e^{ikx}\sin(\alpha_{2}t)\alpha_{2}^{2}(\hat{\Psi}_{2,0}^{(2)}(k) - \hat{\Psi}_{2,0}^{(2)}(-k)) + 2e^{ikx}\sin(\alpha_{2}t)k^{2}(\hat{\Psi}_{2,0}^{(2)}(-k) - \hat{\Psi}_{2,0}^{(2)}(k)) }{4 \alpha_{2} m_{2} \pi i} dk \\
& + \int_{-\infty}^{\infty} \frac{e^{ikx}\cos(\alpha_{2}t)(\alpha_{2}^{2}-k^{2})i(B_{L,2}^{(2)}(i\alpha_{2},t) - B_{L,2}^{(2)}(-i\alpha_{2},t) )(e^{-ikL} - e^{ikL})}{4 \alpha_{2} m_{2} \pi i}dk \\
& + \int_{-\infty}^{\infty} \frac{e^{ikx}\sin(\alpha_{2}t)(\alpha_{2}^{2} - k^{2})(B_{L,2}^{(2)}(i\alpha_{2},t) + B_{L,2}^{(2)}(-i\alpha_{2},t) )(e^{-ikL} - e^{ikL})}{4 \alpha_{2} m_{2} \pi i}dk \\
& + \int_{-\infty}^{\infty} \frac{2e^{ikx}\cos(\alpha_{2}t)m_{2}ki(B_{0,1}^{(2)}(i\alpha_{2},t) - B_{0,1}^{(2)}(-i\alpha_{2},t) )}{4 \alpha_{2} m_{2} \pi i}dk \\
& + \int_{-\infty}^{\infty} \frac{2e^{ikx}\sin(\alpha_{2}t)m_{2}k(B_{0,1}^{(2)}(i\alpha_{2},t) + B_{0,1}^{(2)}(-i\alpha_{2},t) )}{4 \alpha_{2} m_{2} \pi i}dk \\
& - \int_{-\infty}^{\infty} \frac{e^{ikx}m_{2}k\sin(\alpha_{2}t)(B_{L,1}^{(2)}(i\alpha_{2},t) + B_{L,1}^{(2)}(-i\alpha_{2},t)) (e^{ikL} + e^{-ikL})}{4 \alpha_{2} m_{2} \pi i} dk \\
& + \int_{-\infty}^{\infty} \frac{e^{ikx}m_{2}k\sin(\alpha_{2}t)(2\hat{\Psi}_{1,0}^{(2)}(k) + 2\hat{\Psi}_{1,0}^{(2)}(-k))}{4 \alpha_{2} m_{2} \pi i} dk,
\label{Psi_2_1}
\end{split}
\end{eqnarray}
and 
\begin{eqnarray}
\begin{split}
& \Psi_{2}^{(2)}(x,t) =  \int_{-\infty}^{\infty} \frac{e^{ikx}\cos(\alpha_{2}t)m_{2}i(\hat{\Psi}_{1,0}^{(2)}(k) - \hat{\Psi}_{1,0}^{(2)}(-k)) }{2 \pi \alpha_{2} i} dk \\
& + \int_{-\infty}^{\infty} \frac{ e^{ikx}\cos(\alpha_{2}t)\alpha_{2}i(\hat{\Psi}_{2,0}^{(2)}(-k) - \hat{\Psi}_{2,0}^{(2)}(k))}{2 \pi \alpha_{2} i} dk \\
& - \int_{-\infty}^{\infty} \frac{e^{ikx} \cos(\alpha_{2}t)ki \hat{\Psi}_{2,0}^{(2)}(k) + e^{ikx} \sin(\alpha_{2}t)k \hat{\Psi}_{2,0}^{(2)}(-k)}{2 \pi \alpha_{2} i} dk \\
& + \int_{-\infty}^{\infty} \frac{e^{ikx} \cos(\alpha_{2}t) (m_{2} B_{L,1}^{(2)}(-i\alpha_{2},t) + \alpha_{2} B_{L,2}^{(2)}(-i\alpha_{2},t))i(e^{ikL} - e^{-ikL})}{2 \pi \alpha_{2} i} dk \\
& - \int_{-\infty}^{\infty} \frac{e^{ikx} \cos(\alpha_{2}t) k i B_{L,2}^{(2)}(-i\alpha_{2},t)(e^{ikL} + e^{-ikL} - 2)}{2 \pi \alpha_{2} i} dk,
\label{Psi_2_2}
\end{split}
\end{eqnarray}
for $x>0$, $t>0$, and $\alpha_{2}=\sqrt{k^{2}+m_{2}^{2}}$.
\end{proposition}

\section{Conclusion}

\quad In this work, we establish propositions for the solution of the interface problem for the Dirac equation on two semi-infinite domains and two finite domains in the massless and massive cases, respectively. In \cite{11}, we find explicit solutions for the Dirac equation in the massless and massive cases on the positive half line and a finite non-negative interval. The analytical solutions found for interface problem for the Dirac equation in the massless and massive cases that we have and in particular the solution that we obtain for the interface problem for the Dirac equation for the right-hand side (positive half line) matches with the solution for Dirichlet's problem for the Dirac equation on positive half line and finite non-negative interval that appears in \cite{11}, which implies that through the UTM for vector case, we can obtain explicit solutions for interface problems for hyperbolic-type PDEs such as the Klein-Gordon equation (see \cite{Klein_Gordon}) and the Dirac equation. Furthermore, the solutions for an interface problem on two semi-infinite domains and two finite domains in the massive case are convergent explicit integral representations.

\end{document}